\documentclass[times]{nmeauthNoIJNME}

\usepackage{amsmath}
\usepackage{amssymb}
\usepackage{mathtools}
\usepackage[colorlinks,
            bookmarksopen,
            bookmarksnumbered,
            citecolor=blue,
            linkcolor=blue,
            urlcolor=black]{hyperref}

\usepackage{float}
\usepackage{subfig}

\usepackage{csml}

\begin{document}

\runningheads{X Xiao and F Cirak}{Infill topology and shape optimisation of lattice-skin structures}

\title{Infill topology and shape optimisation of lattice-skin structures}

\author{Xiao Xiao\affil{1,2}\corrauth and Fehmi Cirak\affil{1}}

\address{\affilnum{1}Department of Engineering, University of Cambridge, Trumpington Street, Cambridge CB2 1PZ, UK \break
\affilnum{2}Inria, 2004 route des Lucioles, 06902 Sophia Antipolis, France
}

\corraddr{Inria, 2004 route des Lucioles, 06902 Sophia Antipolis, France. E-mail: xiao.xiao@inria.fr}

\begin{abstract}
Lattice-skin structures composed of a thin-shell skin and a lattice infill are widespread in nature and large-scale engineering due to their efficiency and exceptional mechanical properties.  Recent advances in additive manufacturing, or 3D printing, make it possible to create lattice-skin structures of almost any size with arbitrary shape and geometric complexity. We propose a novel gradient-based approach to optimising both the shape and infill of lattice-skin structures to improve their efficiency further. The respective gradients are computed by fully considering the lattice-skin coupling while the lattice topology and shape optimisation problems are solved in a sequential manner. The shell is modelled as a Kirchhoff-Love shell and analysed using isogeometric subdivision surfaces, whereas the lattice is modelled as a pin-jointed truss. The lattice consists of many cells, possibly of different sizes, with each containing a small number of struts. We propose a penalisation approach akin to the SIMP (solid isotropic material with penalisation)  method for topology optimisation of the lattice. Furthermore, a corresponding sensitivity filter and a lattice extraction technique are introduced to ensure the stability of the optimisation process and to eliminate scattered struts of small cross-sectional areas. The developed topology optimisation technique is suitable for non-periodic, non-uniform lattices. For shape optimisation of both the shell and the lattice, the geometry of the lattice-skin structure is parameterised using the free-form deformation technique. The topology and shape optimisation problems are solved in an iterative, sequential manner. The effectiveness of the proposed approach and the influence of different algorithmic parameters are demonstrated with several numerical examples. 
\end{abstract}


\maketitle



%
\section{Introduction \label{sec:intro}}
\begin{figure}
\centering
\subfloat[Geometry (parts of the skin omitted for visualisation).]
{
  \includegraphics[scale=0.95]{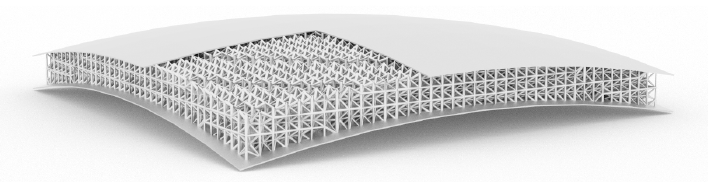}
}
\hfill
\subfloat[Finite element discretisation.]
{
  \includegraphics[scale=0.95]{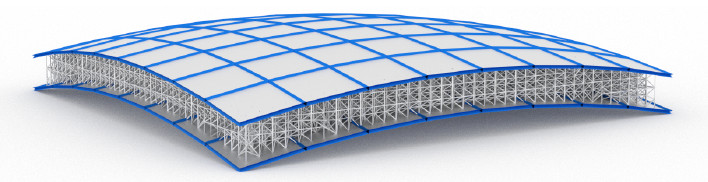}
}
\\
\subfloat[Infill topology and shape optimised structure.]
{
  \includegraphics[scale=0.95]{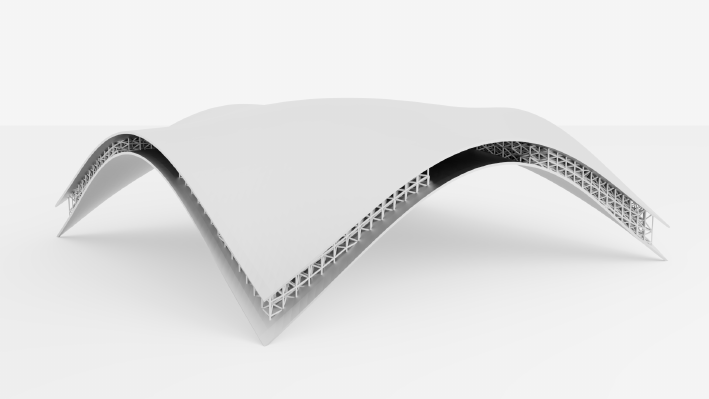}
}
\hfill
\subfloat[Deflected structure (colours represent strut displacements).]
{
  \includegraphics[scale=0.95]{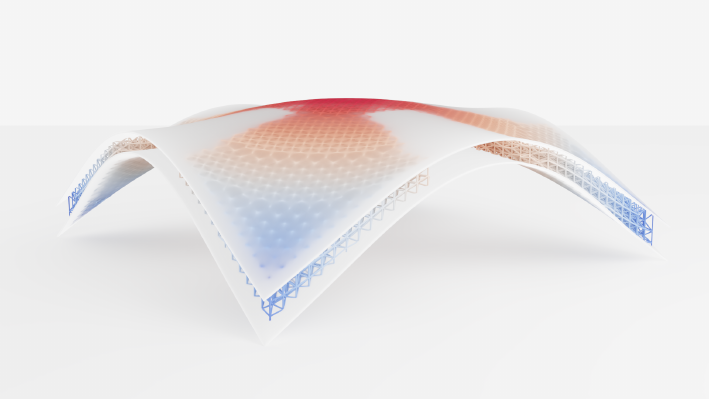}
}
\caption{An illustrative application of the proposed infill topology and shape optimisation approach. The doubly-curved lattice-skin structure (a) with two thin-shell face-sheets and a body-centred cubic (BCC)  lattice core  is subjected to uniform pressure loading. The thin-shells are discretised using isogeometric spline basis functions and the lattice is modelled as a pin-jointed truss (b). The shape and the infill topology of the structure are optimised, in turn, so that its compliance is minimised (c, d).}
\label{fig:latticeSkin}
\end{figure}

%
Composite lattice-skin structures combine the advantages of thin-shell and lattice structures. Owing to their curved geometry, shell structures exhibit a superior load carrying capacity and are lightweight, but respond extremely sensitively to any changes in the loading, geometry, etc.~\cite{bushnell1981buckling,ramm2004shell}. In contrast, lattices are less sensitive, have more easily tunable mechanical properties, and are ideal for multi-functional designs that simultaneously draw on several properties, like stiffness, heat or inertia~\cite{gibson1999cellular,fleck2010micro}. Due to their favourable mechanical properties, lattice-skin structures are prevalent in engineering, e.g. lattice-core sandwich plates (Figure~\ref{fig:latticeSkin}),  and in nature, e.g. trabecular bone. Lattice and lattice-skin structures are currently experiencing a renaissance because of the rapid development of additive manufacturing, or 3D printing, technologies~\cite{gibson2014additive,schaedler2016architected}. Additive manufacturing makes it possible to create parts with nearly any shape and geometric complexity, hierarchically extending from the part-scale down to the resolution length of the printing process. To explore and exploit the resulting vast design space, computational approaches for rational design and optimisation of lattice-skin structures are indispensable. 

Several powerful optimisation approaches have been proposed over the last 50 years for structures consisting only of lattices. Consequently, there is an extensive amount of literature that is impossible to review in detail. The ground structure method is one of the earliest techniques in the topology optimisation of pin-jointed truss structures~\cite{kirsch1989optimal,bendsoe1994optimization, ohsaki2016optimization,zegard2014grand}. It finds an optimal structure by optimising the cross-sectional areas of struts of a given highly-connected truss by usually solving a linear programming problem. The objective is to minimise the total material volume while satisfying the equilibrium and prescribed stress limits. Although in classical formulations the joint positions are fixed, some approaches consider the joint positions as additional optimisation degrees of freedom~\cite{achtziger2007simultaneous}. Broadly, the ground structure method is not scalable for very large problems, even though there are methods to alleviate this limitation by incrementally adding members~\cite{gilbert2003layout, martinez2007growth, hagishita2009topology}. Beyond that, the ground structure method can lead to non-manufacturable thin struts and, crucially, it is unsuitable for optimisation of lattices connected to other structures. Classical homogenisation-based approaches for topology optimisation offer an alternative approach for optimisation of lattices consisting of repetitive unit cells~\cite{bendsoe1988generating, wang2017concurrent, panesar2018strategies}. In such methods, the lattice is first homogenised with the resulting effective material parameters depending on unit cell parameters, like the material content or cell orientation. These parameters are later determined by optimising the continuum structure composed of the homogenised material. Projecting back the obtained parameters to the unit cells gives an optimised lattice structure. In general, adjacent unit cells will not match across common cell boundaries, which has to be ensured by other means~\cite{wang2017concurrent, pantz2008post, groen2018homogenization}. A more fundamental limitation of the homogenisation-based approaches is that the underlying assumption of scale separation breaks down, for instance, near boundaries and when the unit cell size is not sufficiently small compared to the overall size of the structure.  The application of the two mentioned principal lattice optimisation techniques in additive manufacturing context has been explored, amongst others, in~\cite{panesar2018strategies, zegard2016bridging, robbins2016efficient}.

The shape optimisation of shell structures is an equally widely investigated topic. One of the key challenges in any kind of shape optimisation is to find a set of geometric design parameters that can parameterise a sufficiently rich set of geometries. To this end, classical geometry parameterisation techniques make use of splines~\cite{ Braibant:1984aa,han2014adaptive}, filtered/smoothed finite element representations~\cite{le2011gradient, bletzinger2014consistent}, or free-form deformations~\cite{imam1982three,sederberg1986free}. Lately, shape optimisation has greatly benefited from the ongoing academic and industrial interests in isogeometric design and analysis. In isogeometric analysis, the same spline representation is used for finite element analysis and geometric modelling~\cite{hughes2005isogeometric}. We use Catmull-Clark subdivision surfaces and the underlying spline basis functions for  isogeometric analysis.  The Catmull-Clark subdivision surfaces are the generalisation of tensor-product B-splines and NURBS to meshes with arbitrary connectivity~\cite{Catmull:1978aa,wei2015truncated,zhang2018subdivision,wei2021tuned}. Given the unique refinability property of splines, a different resolution of the same spline representation can be used for parameterising the geometry in optimisation~\cite{cirak2002integrated,kiendl2014isogeometric,bandara2018isogeometric, hirschler2019embedded}. An additional benefit of smooth spline basis functions in shell analysis is that they are suitable for discretising the Kirchhoff-Love equations~\cite{Cirak:2000aa,cirak2011subdivision,Long:2012aa}. The resulting discretised equations depend only on the displacement of the shell mid-surface, simplifying the computation of sensitivities, or shape gradients. 

The concurrent shape and topology optimisation of lattice-skin structures has been scarcely explored.  There is somewhat related work on homogenised porous structures with a stiff outer coating, or skin,~\cite{clausen2015topology, clausen2017topology, wu2017minimum, wang2018level}, and on lattice infill optimisation without taking into account the coupling between the lattice and the shell~\cite{wang2013cost}.  To close this gap, we propose a shape and topology optimisation approach by combining topology optimisation of the lattice with shape optimisation of the entire structure. Building on our earlier work on isogeometric design and analysis of lattice-skin structures~\cite{xiao2019interrogation}, the lattice is modelled as a pin-jointed truss and the skin as a Kirchhoff-Love shell. The lattice consists of a large number of cells, which in turn consist of a small number of struts connected by pins that do not transfer moments. This approximation is sufficient for the technologically important stretch-dominated lattice structures, see e.g.~\cite{deshpande2001effective}. The shell and the lattice are coupled using Lagrange multipliers. The topology of the lattice is optimised by taking the cross-sectional areas of the struts as design variables. To obtain an optimised structure with only  struts of desired cross-sectional areas,  we choose, as in the SIMP method~\cite{ bendsoe1999material}, scaled cross-sectional areas as design variables. For instance, for a \emph {black and white} design  with struts of the same cross-sectional area, intermediate cross-sectional areas are penalised using a power function. Furthermore, we introduce a cell-wise defined sensitivity filter to avoid scattered struts similar to checkerboards in continuum structures. The sensitivity filter, in combination with a postprocessing step for cell reconstruction, yields a lattice with no mechanisms~\cite{pellegrino1986matrix}.  In shape optimisation, we consider the entire lattice-skin structure and parameterise its shape using the free-form deformation (FFD) technique~\cite{sederberg1986free}. In the FFD approach the structure is immersed in a control prism, which is in turn parameterised with Bernstein basis functions controlling the shape of the structure via the respective control point positions.  We take the mechanical coupling between the lattice and the shell fully into account in computing the sensitivities for shape and topology optimisation. For ease of implementation and efficiency, the topology and shape optimisation steps are performed sequentially. In the presented examples we first optimise the lattice topology and then the shape of both the shell and the lattice. Although we consider only compliance as a cost function, the proposed approach is straightforward to extend to other cost functions. 

The outline of the paper is as follows. In Section~\ref{sec:analysis}, we briefly review the governing equations of lattice-skin structures and sketch their isogeometric finite element discretisation. Subsequently, in Section~\ref{sec:latticeOpt}, we introduce the proposed lattice topology optimisation technique taking the coupling between the lattice and the skin into account. Building on that, in Section~\ref{sec:shapeOpt}, after briefly introducing the FFD technique for shape parameterisation, we discuss the shape optimisation of lattice-skin structures. We then introduce in Section~\ref{sec:examples} four examples of increasing complexity, with the first three involving the topology optimisation of lattice and lattice-skin structures and the last example combining topology and shape optimisation. 

%
\section{Governing equations and finite element discretisation \label{sec:analysis}}
%
In this section, we briefly review the governing equations of lattice-skin structures and their finite element discretisation, for details refer to~\cite{Cirak:2000aa, xiao2019interrogation}. It is assumed that displacements are small and the material behaviour is elastic. We model the shell as a Kirchhoff-Love thin-shell, the lattice as a pin-jointed truss structure and the two are coupled using Lagrange multipliers.  

The displacement of the shell mid-surface~$\Omega^\mathrm{s}$ is denoted with \mbox{$\vec u^\mathrm{s}  \colon  \Omega^\mathrm{s} \rightarrow \mathbb R^3$} and the displacement of a lattice node, or joint, with the index~$j$ with \mbox{$\vec u_j^\mathrm{l} \in \mathbb R^3$}. The displacement of all lattice nodes is collected in the matrix~$\vec U^\mathrm{l} \in \mathbb{R}^{n^\mathrm{l} \times 3}$, where~$n^\mathrm{l}$ is the total number of nodes. Henceforth, the superscripts $\mathrm{s}$ and $\mathrm{l}$ are used to distinguish between the thin-shell and the lattice variables. The total potential energy~$\Pi(\vec{u}^\mathrm{s}, \vec{U}^\mathrm{l})$ of the lattice-skin structure is composed of the thin-shell energy~$\Pi^\mathrm{s}(\vec{u}^\mathrm{s})$, the lattice energy~$\Pi^\mathrm{l}(\vec{U}^\mathrm{l})$, and the potential of the applied loads~$\Pi^\mathrm{ext}(\vec{u}^\mathrm{s})$, i.e.
\begin{equation} \label{eq:energyFunctional} 
	\Pi(\vec{u}^\mathrm{s}, \vec{U}^\mathrm{l}) = \Pi^\mathrm{s}(\vec{u}^\mathrm{s}) + \Pi^\mathrm{l}(\vec{U}^\mathrm{l}) + \Pi^\mathrm{ext}(\vec{u}^\mathrm{s})  \, . 
\end{equation}

The lattice energy is comprised of the internal energies of the individual struts,
\begin{equation}
	\Pi^\mathrm{l}(\vec{U}^\mathrm{l}) = \sum_{e}W^{\mathrm l} \left ({\epsilon}_e \right )A_e l_e \, , 	
\end{equation}
where~$W^{\mathrm l}({\epsilon}_e)$ is the internal energy density depending on the axial strain~${\epsilon}_e (\vec{U}^\mathrm{l})$,~$A_e$ is the cross-sectional area, and~$l_e$ is the length of the strut with the index~$e$.  It is assumed that the cross-sectional area of the struts is constant over their lengths. The internal energy density~$W^{\mathrm l}({\epsilon}_e)$ for an elastic material with Young's modulus~$E$ is given by 
\begin{equation}
W^{\mathrm l}(\epsilon_e) = \frac{1}{2}E\epsilon_e^2 \, .
\end{equation}

The potential energy of the displaced thin-shell takes the form
\begin{equation}
	\Pi^\mathrm{s}(\vec{u}^\mathrm{s}) = \int_{\Omega^\mathrm{s}}\! \left (W^{\mathrm{s,m}}(\vec{\alpha}) + W^{\mathrm{s,b}}(\vec{\beta}) \right ) \,\D \Omega^\mathrm{s}  \, , 
\end{equation}
where~$W^{\mathrm{s,m}}(\vec{\alpha})$ and~$W^\mathrm{s,b}(\vec{\beta})$ are the membrane and bending strain energy densities depending on the membrane  and bending strain tensors~$\vec{\alpha}(\vec u^\mathrm{s})$ and $\vec{\beta} (\vec u^\mathrm{s}) $, respectively. For a thin-shell with the thickness~$t$ and an elastic and isotropic material with~Young's modulus~$E$ and Poisson's ratio~$\nu$, the two energy densities are given by 
\begin{equation}
W^{\mathrm{s,m}}(\vec{\alpha}) = \frac{1}{2}\frac{Et}{1 -
\nu^2}\vec{\alpha}:\vec{H}:\vec{\alpha} \, , \quad
W^{\mathrm{s,b}}(\vec{\beta}) = \frac{1}{2}\frac{Et^3}{12(1 -
\nu^2)}\vec{\beta}:\vec{H}:\vec{\beta} \, , 
\end{equation}
where $\vec{H}$ is a fourth-order geometry-dependent tensor~\cite{Cirak:2000aa}. Finally, assuming, without loss of generality, that the only applied external loading is a distributed load~$\vec p \colon  \Omega^\mathrm{s} \rightarrow \mathbb R^3$, the potential of the applied loads takes the form 
\begin{equation}
	\Pi^\mathrm{ext}(\vec{u}^\mathrm{s}) = 	- \int_{\Omega^\mathrm{s}} \vec p \cdot \vec u^\mathrm{s} \, \D \Omega^\mathrm{s} \, . 
\end{equation}
\begin{figure}
\centering
\includegraphics[scale = 0.8]{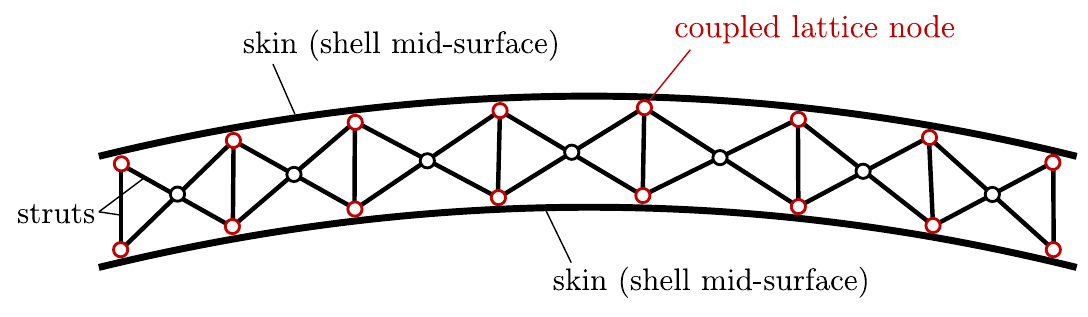}
\caption{Schematic of a lattice-skin structure consisting of two thin-shells and a lattice core. Some of the lattice nodes (red circles) are constrained to have the same displacement as the shell mid-surface. The small offset between the joint position and its projection on the mid-surface is neglected because of the small shell thickness. }
\label{fig:latticeSkinAnalysis}
\end{figure}

The displacements of the lattice nodes attached to the shell are required to be compatible with the shell displacement, see Figure~\ref{fig:latticeSkinAnalysis}. To enforce it, first  we collect the indices of all lattice nodes attached to the shell in the set $\mathcal{D}_\mathrm{c}^\mathrm{l}$. Furthermore, as usual in geometrically exact shell theories~\cite{simo1989stress,ciarlet2005introduction}, the shell mid-surface is  parameterised with curvilinear parametric coordinates~$\vec \theta \in \mathbb R^2$. Hence, the compatibility of the lattice node displacements and the shell displacement requires  
\begin{equation} \label{eq:dispCompatibility}
\vec{u}^\mathrm{s}(\vec{\theta}_j^\mathrm{l}) =
\vec{u}^\mathrm{l}_j \quad \forall j\in\mathcal{D}_\mathrm{c}^\mathrm{l} \, ,
\end{equation}
where $\vec{u}_j^\mathrm{l}$ is the displacement of the $j$-th lattice node attached to the shell mid-surface at the point with the parametric coordinate~$\vec{\theta}_j^\mathrm{l}$.  This compatibility constraint~\eqref{eq:dispCompatibility} is imposed with Lagrange multipliers $\vec{\lambda}_j$ via the augmented potential
\begin{equation} \label{eq:lagrangian}
L(\vec{u}^\mathrm{s}, \vec{U}^\mathrm{l}, \vec{\lambda}_j) = \Pi(\vec{u}^\mathrm{s},
\vec{U}^\mathrm{l}) + \sum_{j \in \set{D}_\mathrm{c}^
\mathrm{l}}\vec{\lambda}_{j}\left( \vec{u}^\mathrm{l}_j - \vec{u}^\mathrm{s}(\vec{\theta}_j^\mathrm{l}) \right) \, .
\end{equation}
The stationarity of this potential yields the following set of equilibrium and compatibility equations
\begin{subequations}
\begin{align}
\frac{\partial\Pi^{\mathrm{s}}(\vec{u}^\mathrm{s})}{\partial\vec{u}^\mathrm{s}}\delta\vec{u}^\mathrm{s}
+
\frac{\partial\Pi^{\mathrm{ext}}  (\vec{u}^\mathrm{s}) }{\partial\vec{u}^\mathrm{s}}\delta\vec{u}^\mathrm{s}
-
\sum_{j \in \mathcal{D}_\mathrm{c}^\mathrm{l}} \vec{\lambda}_j\delta\vec{u}^\mathrm{s}(\vec{\theta}_j^\mathrm{l}) &= 0 \, , \label{eq:varshell}\\[2pt]
\frac{\partial\Pi^{\mathrm{l}}(\vec{U}^\mathrm{l})}{\partial\vec{u}_j^\mathrm{l}}  &= 0 \quad
\forall j \notin\mathcal{D}_\mathrm{c}^\mathrm{l} \, , \\[2pt]
\frac{\partial\Pi^{\mathrm{l}}(\vec{U}^\mathrm{l})}{\partial\vec{u}_j^\mathrm{l}}
+ \vec{\lambda}_j &= 0 \quad \forall j \in \mathcal{D}_\mathrm{c}^\mathrm{l} \, ,
\label{eq:varlattice}\\[2pt]
\vec{u}^\mathrm{l}_j - \vec{u}^\mathrm{s}(\vec{\theta}_j^\mathrm{l})  &= 0 \quad \forall j \in\mathcal{D}_\mathrm{c}^\mathrm{l} \, , 
\end{align}
\label{eq:varLagrangian}
\end{subequations}
where~$\delta \vec u^{\mathrm s}$ denotes the virtual displacement of the shell mid-surface. 

We discretise the mid-surface of the thin-shell~$ \Omega^\mathrm{s}$  and its displacement~$\vec u^{\mathrm s}$ with Catmull-Clark subdivision basis functions, which represent a generalisation of cubic tensor-product B-splines to unstructured meshes; see~\cite{Cirak:2000aa,zhang2018subdivision} for details.  After introducing the discretised mid-surface and the discretised displacements into the weak form~\eqref{eq:varLagrangian}, and subsequent numerical evaluation of integrals, we obtain the  discrete system of equations 
\begin{equation} \label{eq:equilibrium}
\begin{pmatrix}
\,\,\,\ary{K}^\mathrm{s} & \ary{0} & -\ary{G}_\mathrm{c}^\trans \\[2pt]
\,\ary{0} & \,\ary{K}^\mathrm{l} & \,\,\,\ary{I}_\mathrm{c}^\trans \\[2pt]
-\ary{G}_\mathrm{c} & \,\ary{I}_\mathrm{c} & \ary{0}
\end{pmatrix}
\begin{pmatrix}
\ary{u}^\mathrm{s} \\[2pt]
\ary{u}^\mathrm{l} \\[2pt]
\ary{\lambda}
\end{pmatrix}
=
\begin{pmatrix}
\ary{f}^\mathrm{s} \\[2pt]
\ary 0 \\[2pt]
\ary{0}
\end{pmatrix} \, ,
\end{equation}
where $\ary{K}^\mathrm{s}$ and $\ary{K}^\mathrm{l}$ are the stiffness matrices of the thin-shell and the lattice; $\ary{u}^\mathrm{s}$ and $\ary{u}^\mathrm{l}$ are the respective displacement vectors and $\ary{f}^\mathrm{s}$ is the external force vector of the shell; $\ary{I}_\mathrm{c}$ is an extraction matrix which multiplied with $\ary{u}_\mathrm{c}^\mathrm{l}$ gives the displacements of the lattice nodes with $j \in \set{D}_\mathrm{c}^\mathrm{l}$, and $\ary{G}_\mathrm{c}$ is a matrix obtained by evaluating the subdivision basis functions at the parametric mid-surface coordinates~$\vec \theta^{\mathrm l}_j$ corresponding to lattice nodes with~$j\in\mathcal{D}_\mathrm{c}^\mathrm{l}$. We write~\eqref{eq:equilibrium} more compactly as 
\begin{equation}  \label{eq:equilibrium2}
	\ary K \ary u = \ary f \, . 
\end{equation}

%
\section{Lattice topology optimisation \label{sec:latticeOpt}}
%

%
\subsection{Problem statement and penalisation of cross-sectional areas}
%
We optimise the topology of the lattice structure with a penalisation technique akin to the SIMP approach extensively used in topology optimisation of continuum structures~\cite{bendsoe1999material}. In this section we focus on optimised black and white designs with struts of the same cross-sectional areas. If an optimised design with cross-sectional areas within a certain range is desired, the presented approach can be extended as demonstrated in Section~\ref{sec:cantileverLatticeBezier}.
In a black and white design intermediate cross-sectional areas are penalised to obtain an optimised lattice with only some of the struts removed and (almost) no changes to the cross-sectional areas of the remaining ones. The associated compliance topology optimisation problem considering the equilibrium equation~\eqref{eq:equilibrium} reads
\begin{subequations} \label{eq:topOpt}
\begin{align}
\displaystyle\underset{\ary \rho}{\text{minimise}} \quad  & J( \ary \rho) = \ary f^\mathrm{s} \cdot \ary u^\mathrm{s} = \ary u^\mathrm{s} (\ary \rho) \cdot \ary K^{\mathrm s} \ary u^\mathrm{s} (\ary \rho) + \ary u^\mathrm{l} (\ary \rho) \cdot \ary K^{\mathrm l} (\ary \rho) \ary u^\mathrm{l} (\ary \rho) \, , \label{eq:topOptA}\\
\text{subject to} \quad  & \ary K (\ary \rho) \ary u  = \ary f \, , \label{eq:topOptB}\\
& V^\mathrm{l} / \, \overline V^\mathrm{l} \le V_f^\mathrm{l} \, , \label{eq:topOptC} \\
& \ary 0 \le \ary \rho \le  \ary 1 \, ,  \label{eq:topOptD}
\end{align}
\end{subequations}
where $J(\ary{\rho})$ is the structural compliance of the lattice-skin structure, $\ary \rho$ is the relative density vector yet to be defined; $V^\mathrm{l}$ is the actual and $\overline V^\mathrm{l}$ the initial lattice material volume, and $V_f^\mathrm{l}$ is the prescribed volume fraction of the lattice. Note that in~\eqref{eq:topOptA}, according to~\eqref{eq:equilibrium}, only the shell has an external loading and the shell stiffness matrix $\vec K^\mathrm{s}$ is independent of the relative density~$\vec \rho$. 

\begin{figure}
\centering
  \includegraphics[scale=1.4]{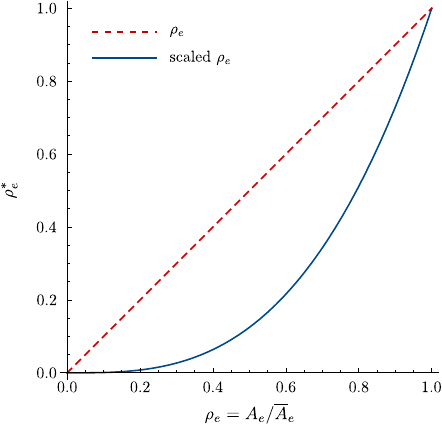}
  \caption{The penalised relative density$\rho_e^* = \rho_e^p$ with $p = 3$  and the area constraint $A_e \leq \overline{A}_e$.}   \label{fig:simpExp}
\end{figure}

%
%

In our SIMP-like penalisation approach the element stiffness matrix of each strut is first expressed as 
\begin{equation} \label{eq:areaDensity}
\ary{K}_e^\mathrm{l} = \frac{A_e}{\overline{A}_e}\overline{\ary{K}}_e^\mathrm{l} = \rho_e\overline{\ary{K}}_{e}^\mathrm{l} \, ,
\end{equation}
where $\overline{\ary{K}}_{e}^\mathrm{l}$ is the original stiffness matrix of the strut with the index~$e$, and its ratio of the current and original cross-sectional areas $A_{e}$ and $\overline{A}_{e}$ is defined as its  relative density~$\rho_e$. Hence, the relative density may take the values~\mbox{$ 0 \le \rho_e \le 1$}, in line with~\eqref{eq:topOptD}. To penalise the intermediate relative densities~\mbox{$ 0 < \rho_e  <1$} the strut element stiffness matrix is replaced with
\begin{equation} \label{eq:stiffnessPenalty}
\ary{K}_e^\mathrm{l}(\rho_e^*) = \rho_e^*\overline{\ary{K}}_{e}^\mathrm{l} \, ,
\end{equation}
with the penalised relative density
\begin{equation} \label{eq:scaling}
\rho_e^* = 
\begin{cases}
\; \rho_e^p \, &\text{if} \; \rho_e < 1 \, , \\
\; \rho_e \, &\text{otherwise} \, .
\end{cases}
\end{equation}
The penalisation parameter~$ p \gtrsim 3$ ensures that the struts with relative densities close to~\mbox{$\rho_e = 0$} and~\mbox{$\rho_e=1$} are preferred; see Figure~\ref{fig:simpExp}. In some applications it can be advantageous or desirable to have relative densities~$\rho_e >1$ so that the optimised area~\mbox{$A_e$ can be larger than the original area $ \overline{A}_{e}$.} In such cases the definition~\eqref{eq:scaling} leads to a slope discontinuity at~$\rho_e = 1$, which may affect the convergence of optimisation. To remedy this, the power function~$\rho_e^p$ can be replaced, for instance, with a B\'ezier curve with continuous derivatives at~$\rho_e = 1$. 
%
 In the examples presented in this paper, we consider only the case~$\rho_e \le 1$ and apply this as a constraint during optimisation. 
 In addition, we always apply the constraint~\mbox{$\rho_e > \rho_{\mathrm{min}} \approx10^{-6} $} to avoid the singularity of the global stiffness matrix. As discussed in Section~\ref{sec:cantileverLatticeBezier}, the power function~$\rho_e^p$ in~\eqref{eq:scaling}  can be replaced with a suitably shaped B\'ezier spline curve  to obtain  cross-sectional areas within a certain desired range.

%
\subsection{Sensitivity analysis and filtering}
%
For gradient-based optimisation the derivatives, or sensitivities, of the cost function and the constraints in~\eqref{eq:topOpt} with respect to the cross-sectional areas are required. It is straightforward to show that the derivative of the cost function~\eqref{eq:topOptA} with respect to the penalised densities, taking into account the equilibrium equation~\eqref{eq:equilibrium}, is given by 
\begin{equation}
\frac{\partial J(\ary{\rho}^* )}{\partial \rho_e^*} = -\ary{u}_e^\mathrm{l}\cdot\frac{\partial
\ary{K}_e^\mathrm{l}(\rho_e^*)}{\partial \rho_e^*}\ary{u}_e^\mathrm{l} \, .
\end{equation}
After introducing the definitions of the relative density and the penalised stiffness matrix, i.e. \eqref{eq:areaDensity} and \eqref{eq:stiffnessPenalty}, we obtain
\begin{equation} \label{eq:latticeDerivative}
\frac{\partial J(\ary{\rho}^*)}{\partial A_e} =
-\frac{\partial \rho_e^*}{\partial
\rho_e}\ary{u}_e^\mathrm{l}\cdot\frac{\overline{\ary{K}}_{e}^\mathrm{l}}{\overline{A}_{e}}\ary{u}_e^\mathrm{l} \, .
\end{equation}

\begin{figure}
\centering
\subfloat[Lattice.]
{
  \includegraphics[scale=0.9]{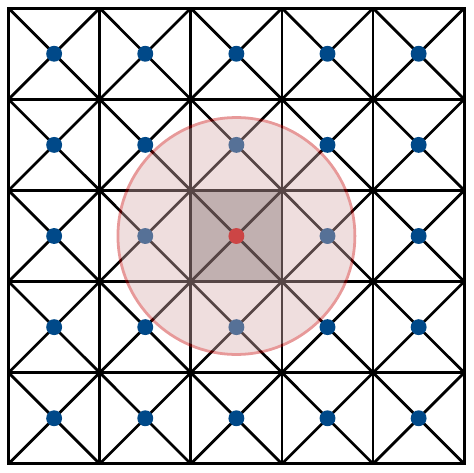}
  \label{fig:filtering1}
}
\hspace{1.5cm}
\subfloat[Unit cells within filter support.]
{
  \qquad\includegraphics[scale=0.9]{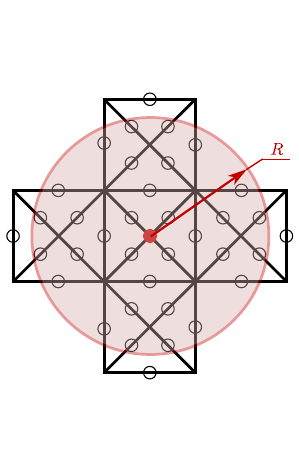}\qquad
  \label{fig:filtering2}
}
\caption{Computation of the filtered sensitivity of a unit cell. The centroid of the considered unit cell is denoted with the red dot. In (a) the solid dots are the centroids of unit cells and the shaded disk is the support of the filter. In (b) the empty dots represent centroids of struts. Only five unit cells are considered within the support of the filter shown. }
\label{fig:sensitivityFilter}
\end{figure}

As known in topology optimisation of continuum structures, filtering techniques are needed to avoid checkerboard instabilities and excessive mesh dependency of the solution~\cite{sigmund1998numerical, sigmund2007morphology}. Similar issues can be observed in lattice topology optimisation without filtering. Commonly, filtering is applied by convolving the computed sensitivities with a kernel, or filter.  We  propose a sensitivity filter for lattices consisting of unit cells, as illustrated in Figure~\ref{fig:sensitivityFilter}. The support of the filter is defined by a prescribed filter radius~$R$ describing a  circle (in 2D) or a sphere (in 3D).  A unit cell is considered within the filter support when its centroid lies within the support, see Figure~\ref{fig:filtering1}. The filtered sensitivity of a unit cell is obtained by centering the filter at its centroid and calculating the weighted average of  the sensitivities of the struts belonging to unit cells within the support. Hence, the filtered sensitivity of a unit cell~$c$ is given by   
\begin{equation} \label{eq:sensitivityFilter}
\widehat{\frac{\partial J(\ary \rho^*)}{\partial A_c}} = \dfrac{\displaystyle\sum_{e=1}^{n_c}
w_e\dfrac{\partial J}{\partial A_e}\dfrac{1}{l_e}}{\displaystyle\sum_{e=1}^{n_c} \dfrac{w_e}{l_e}} \, ,
\end{equation}
where $n_c$ is the number of struts in the unit cells within the support, and $w_e$ is a weight according to a linearly decaying kernel function 
\begin{equation}
w_e(\vec{x}_e, \vec{x}_c) = \max\left\{R - \dist(\vec{x}_e,
\vec{x}_c), \, 0\right\} \, ,
\end{equation}
where $\vec{x}_e$ and $\vec{x}_c$ are the coordinates of the centroid of the $e$-th strut and the cell~$c$, respectively. The support of the filter should not be chosen smaller than a unit cell. The motivation for the division by the strut length~$l_e$ in~\eqref{eq:sensitivityFilter} is that the summation in the numerator can thus be interpreted as the weighted strain energy density, i.e. the strain energy per unit volume. It is worth emphasising that the proposed  filter is based on unit cell sensitivities rather than directly on strut sensitivities. As also reported in~\cite{xia2013method}, directly using the strut sensitivities without filtering easily leads to non-grid lattice layouts and mechanisms.
%
\subsection{Lattice extraction and reconstruction}
\begin{figure}
\centering
\subfloat[Struts with a sensitivity above a threshold. \label{fig:latticeExtractionA}]
{
  \includegraphics{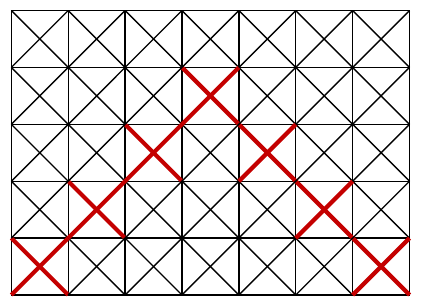}
}
\hfill
\subfloat[Extracted struts with the  recovered unit cells.  \label{fig:latticeExtractionB}]
{
  \includegraphics{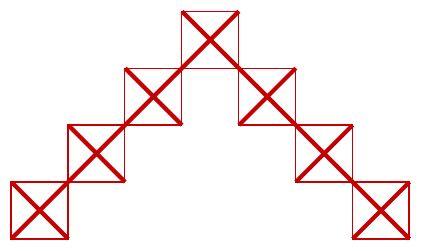}
}
\hfill
\subfloat[Additional struts in the concave void regions. \label{fig:latticeExtractionC}]
{
  \includegraphics{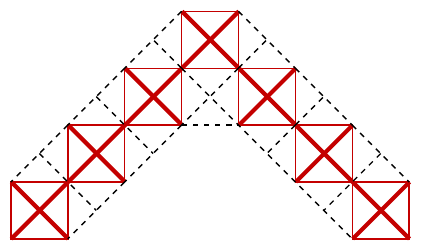}
}
\caption{Extraction and reconstruction of a lattice from topology optimisation result.}
\label{fig:latticeExtraction}
\end{figure}

%
We assign the filtered sensitivity of a unit cell to each strut within the cell. The sensitivity of the struts belonging to several unit cells is obtained by averaging the respective unit cell sensitivities. To obtain the optimised lattice structure, struts with cross-sectional areas larger than a small user-defined threshold are extracted from the lattice, see Figure~\ref{fig:latticeExtractionA}. Notice that struts unique to a unit cell are either all present or not in the optimised topology because they have the same sensitivities.
%
Subsequently, as illustrated in Figure~\ref{fig:latticeExtractionB} the complete topology of a unit cell is recovered when there are any  dangling diagonal struts in the optimised lattice. As the last step, additional struts are reintroduced in the concave void regions of the lattice, see Figure~\ref{fig:latticeExtractionC}. The proposed extraction technique yields a rigid lattice with no mechanisms,  or zero energy modes~\cite{pellegrino1986matrix, deshpande2001foam}. In this work we consider only unit cells that are not mechanisms, or do not have zero energy modes, when considered individually.

%
\section{Shape optimisation of the lattice-skin structure \label{sec:shapeOpt}}
%
\subsection{Geometry parameterisation}
%
\begin{figure}
\centering
\includegraphics[width=0.7\textwidth]{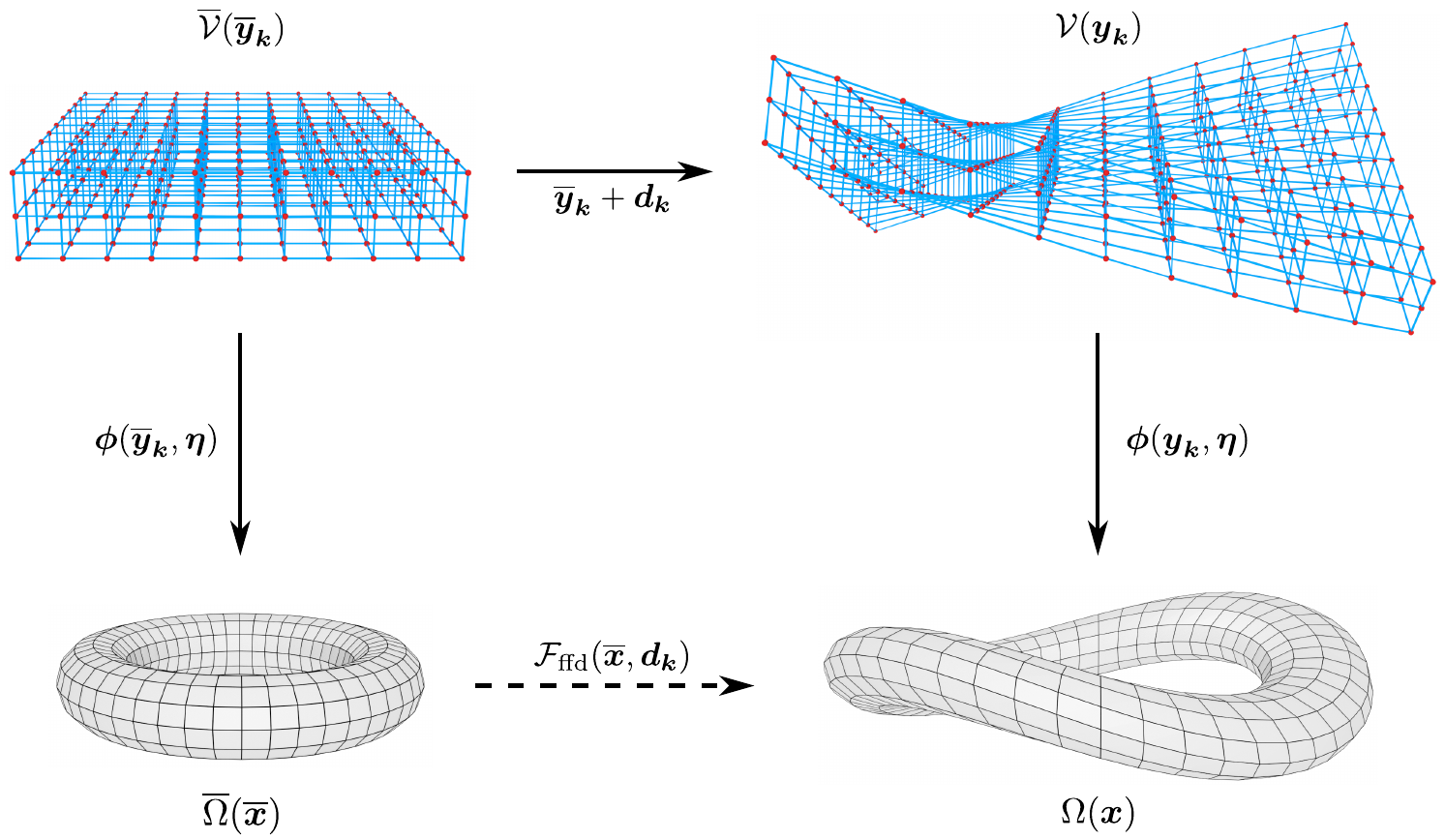}
\caption{Maps involved in shape optimisation using the free-form deformation technique. The geometry of the torus structure  (bottom left) is altered by updating the control point coordinates of the control prism (top left).}
\label{fig:ffd}
\end{figure}

We use the free-form deformation (FFD) technique~\cite{sederberg1986free} to parameterise the overall shape of the lattice-skin structure. The structure with the physical domain $\overline \Omega \subset \mathbb R^3$ is first immersed in a larger rectangular control prism~\mbox{$\mathcal{\overline V}  \supset \overline \Omega $}, see Figure~\ref{fig:ffd}.  An overbar denotes, here and in the following, domains and points related to the original, i.e. not optimised, lattice-skin structure. The control prism is discretised with a uniform grid consisting of \mbox{$[\mu_1 + 1] \times [\mu_2 + 1] \times [\mu_3 + 1]$} control points, where~$\mu_1, \, \mu_2,  \, \mu_3 \ge 1$. The coordinates of the control points are denoted with~$\vec{\overline y}_{\vec{k}} \in \mathbb R^3$ with the multi-index  \mbox{$\vec{k} = (k_1,\,k_2,\,k_3)$}. Subsequently, the control prism is parameterised using trivariate Bernstein basis functions~\mbox{$B^{\vec \mu}_{\vec k} (\vec \eta)$} of degree \mbox{$\vec \mu = (\mu_1,\,\mu_2,\,\mu_3)$}. For the sake of simplicity in presentation, in the following we choose \mbox{$\mu_1 = \mu_2 = \mu_3 = \mu$}. 

In shape optimisation the control point coordinates serve as design degrees of freedom. Each point with the coordinate~$\vec{ \overline x}$ in the original physical domain~$\overline \Omega$ has a fixed parametric coordinate \mbox{$\vec{\eta} = (\eta^1, \eta^2, \eta^3)$}, which is straightforward to determine given that a uniform grid is used. Hence, the introduced Bernstein basis functions provide the following parameterisation of the physical domain 
\begin{equation} \label{eq:ffdCoords}
\vec{\overline x} = \vec{\phi}( \,\vec{\overline y}_{\vec k}, \vec{\eta}) =
\sum_{\vec{k}}^{\vec{\mu}}B_{\vec{k}}^{\vec{\mu}}(\vec{\eta})  \vec{\overline y }_{\vec{k} }
= \sum_{k_1 = 1}^{\mu + 1}\sum_{k_2 = 1}^{\mu + 1}\sum_{k_3 = 1}^{\mu + 1}B_{k_1}^{\mu}(\eta^1)B_{k_2}^{\mu}(\eta^2)B_{k_3}^{\mu}(\eta^3) \vec{\overline y}_{\vec{k}}
\, . 
\end{equation}
The overall shape of the physical domain~$\overline\Omega$ is altered by displacing the control points~$\vec{\overline y} _{\vec{k}}$. Applying a displacement~$\vec d_{\vec k}$ to the original control points  yields the displaced coordinates~$\vec{y}_{\vec k} = \vec{\overline y}_{\vec k} + \vec d_{\vec k}$, and, in turn, this yields the deformed coordinates~$\vec x \in \Omega$ for points in the physical domain according to~\eqref{eq:ffdCoords}. For the sake of completeness, the univariate Bernstein polynomials in~\eqref{eq:ffdCoords} are given by 
\begin{equation}
B_k^\mu(\xi) = \binom{\mu}{k - 1}\xi^{k - 1}(1 - \xi)^{\mu - k + 1} \, , \quad
\xi\in[0,\,1] \, .
\end{equation}

%
\subsection{Problem statement and sensitivity analysis}
%
In shape optimisation we aim to find the optimal positions of the control points $\vec y_{\vec{k}}$ which minimise the compliance of the lattice-skin structure whilst satisfying all constraints. Hence, the shape optimisation problem can be formulated as 
\begin{subequations} \label{eq:shapeOpt}
\begin{align}
\displaystyle\underset{\vec y_{\vec k}}{\text{minimise }} \quad  & J\left(\vec{x}^\mathrm{s}(\vec y_{\vec{k}}),
\vec{x}^\mathrm{l}(\vec y_{\vec{k}})\right) =  \ary u^\mathrm{s} (\vec y_{\vec{k}}) \cdot \ary K^\mathrm{s} (\vec y_{\vec{k}}) \ary u^\mathrm{s}(\vec y_{\vec{k}}) +  \ary u^\mathrm{l} (\vec y_{\vec{k}}) \cdot \ary K^\mathrm{l} (\vec y_{\vec{k}}) \ary u^\mathrm{l}(\vec y_{\vec{k}}) \, , \\
\text{subject to} \quad  & \ary K (\vec{x}^\mathrm{s}, \vec{x}^\mathrm{l}) \ary u 
= \ary f \, , \\
& V \le \overline{V} \, ,
\end{align}
\end{subequations}
where $V$ and $\overline V$ are the actual and the initial total lattice-skin material volume, respectively. The derivative, or sensitivity, of the compliance cost function~$ J(\vec{x}^\mathrm{s}(\vec y_{\vec{k}}), \vec{x}^\mathrm{l}(\vec y_{\vec{k}}) ) $ with respect to the control point coordinates $\vec y_{\vec k}$, taking into account~\eqref{eq:equilibrium2}, reads
\begin{equation} \label{eq:latticeSkinShapeDerivElement}
\frac{\partial J(\vec{x}^\mathrm{s}, \vec{x}^\mathrm{l})}{\partial
\vec{y}_{\vec{k}}} = -\sum_{e}
\ary{u}_e^\mathrm{s} \cdot \frac{\partial\ary{K}_e^\mathrm{s}}{\partial \vec{y}_{\vec{k}}}\ary{u}_e^\mathrm{s}
-
\sum_{e} \ary{u}_e^\mathrm{l} \cdot \frac{\partial\ary{K}_e^\mathrm{l}}{\partial \vec{y}_{\vec{k}}}\ary{u}_e^\mathrm{l}
\, ,
\end{equation}
where~$\ary{K}_e^\mathrm{s}$ and~$\ary{K}_e^\mathrm{l}$ are the shell and strut element stiffness matrices, and  the summations are over the respective elements in the discretised structure. The derivative of the stiffness matrix of a shell element with the index $e$ is given by
\begin{equation} \label{eq:shellDeriv}
\frac{\partial\ary{K}_e^\mathrm{s}}{\partial \vec{y}_{\vec{k}}} =
\sum_i\frac{\partial \ary{K}_e^\mathrm{s}}{\partial
\vec{x}_i^\mathrm{s}}\frac{\partial
\vec{x}_i^\mathrm{s}(\vec{\eta}_i^\mathrm{s})}{\partial \vec{y}_{\vec{k}}} =
\sum_i\frac{\partial \ary{K}_e^\mathrm{s}}{\partial
\vec{x}_i^\mathrm{s}}B_{\vec{k}}^{\vec{\mu}}(\vec{\eta}_i^\mathrm{s}) \, ,
\end{equation}
where the summation is over the nodes  of the $e$-th shell element, and $\vec \eta_i^\mathrm{s}$ is the parametric control prism coordinate of the node with the index $i$. Note that in subdivision surfaces the number of nodes of an element depends on the local connectivity of the mesh. The derivatives of the strut element stiffness matrices have a similar form and are given by
\begin{equation} \label{eq:strutDeriv}
\frac{\partial\ary{K}_e^\mathrm{l}}{\partial \vec{y}_{\vec{k}}} = \sum_{i=1}^2 \frac{\partial\ary{K}_e^\mathrm{l}}{\partial\vec{x}_i^\mathrm{l}}\cdot\vec x^{\mathrm{l}}_{i, \vec y_{\vec{k}}} \, ,
\end{equation}
where the summation is now over the two nodes of the $e$-th strut. Here, it is necessary to distinguish between the lattice nodes attached to the shell and those not. The derivatives of a lattice node $\vec x^{\mathrm{l}}_j$ with respect to the control points $\vec y_{\vec{k}}$ are calculated as
\begin{equation}
\vec x^{\mathrm{l}}_{j, \vec y_{\vec{k}}} = 
\begin{cases}
\displaystyle
\frac{\partial
\vec{x}(\vec \eta_j^\mathrm{l})}{\partial \vec{y}_{\vec{k}}} = B_{\vec{k}}^{\vec{\mu}}(\vec{\eta}_j^\mathrm{l}) \quad \text{if } j \notin \mathcal{D}_{\mathrm{c}}^{\mathrm{l}} \, , \\\noalign{\vskip9pt}
\displaystyle
\sum_i\frac{\partial
\vec{x}^\mathrm{s}(\vec{\theta}_j^\mathrm{l})}{\partial \vec{x}_i^\mathrm{s}}
\frac{\partial \vec{x}_i^\mathrm{s}}{\partial \vec{y}_{\vec{k}}} = 
\sum_i N_i^\mathrm{s}(\vec{\theta}_j^\mathrm{l})
B_{\vec{k}}^{\vec{\mu}}(\vec{\eta}_i^\mathrm{s})
 \quad \text{if } j \in \mathcal{D}_{\mathrm{c}}^{\mathrm{l}} \, ,
\end{cases}
\end{equation}
where the summation is over the nodes~$\vec x_i^\mathrm{s}$ of the shell element to which the lattice node is attached. Furthermore, it is assumed that the shell mid-surface is interpolated with $\vec x^\mathrm{s} (\vec \theta) = \sum_i N_i^\mathrm{s}(\vec \theta ) \vec x_i^\mathrm{s}$ where~$N_i^\mathrm{s}$ are the subdivision basis functions. In order to enforce the displacement coupling in shape optimisation, $\vec x^{\mathrm{l}}_j \in \mathcal{D}_{\mathrm{c}}^{\mathrm{l}}$ should be evaluated with the parameter $\vec \theta_j^\mathrm{l}$ in the attached shell element, while the positions of other lattice nodes are evaluated with the corresponding parameter $\vec \eta$ in the control prism $\mathcal{V}$. In this manner the conformality between the lattice and the thin-shell is guaranteed to be maintained. For derivatives of the shell and lattice stiffness matrices with respect to nodal coordinates, in~\eqref{eq:shellDeriv} and~\eqref{eq:strutDeriv}, we refer to~\cite{bandara2018isogeometric, yin2020topologically}.

The sensitivity of the total volume $V(\vec x^\mathrm{s}, \vec x^\mathrm{l})$ with respect to control point coordinates $\vec y_{\vec k}$ is composed of 
\begin{equation}
\frac{\partial V(\vec x^\mathrm{s}, \vec x^\mathrm{l})}{\partial \vec y_{\vec k}} = \frac{\partial V^\mathrm{s}(\vec x^\mathrm{s})}{\partial \vec y_{\vec k}} + \frac{\partial V^\mathrm{l}(\vec x^\mathrm{l})}{\partial \vec y_{\vec k}} = \frac{\partial A^\mathrm{s}(\vec x^\mathrm{s})}{\partial \vec y_{\vec k}}t^\mathrm{s} + \sum_e \frac{\partial l_e(\vec x^\mathrm{l})}{\partial \vec y_{\vec k}}A_e \, ,
\end{equation}
with the thin-shell and lattice volumes $V^\mathrm{s}$ and $V^\mathrm{l}$, the mid-surface area~$A^\mathrm{s}$ and the thin-shell thickness~$t^\mathrm{s}$. The derivative of the mid-surface area is given by
\begin{equation}
\frac{\partial A^\mathrm{s}(\vec x^\mathrm{s})}{\partial \vec y_{\vec k}} = \sum_i \frac{\partial A^\mathrm{s}}{\partial \vec x_i^\mathrm{s}}\frac{\partial \vec x_i^\mathrm{s}(\vec \eta_i^\mathrm{s})}{\partial \vec y_{\vec k}} = \sum_i \frac{\partial A^\mathrm{s}}{\partial \vec x_i^\mathrm{s}} B_{\vec k}^{\vec \mu}(\vec \eta_i^\mathrm{s}) \, ,
\end{equation}
and the derivative of the strut lengths by
\begin{equation}
\frac{\partial l_e}{\partial \vec y_{\vec k}} = \sum_{i =1}^2 \frac{\partial l_e}{\partial \vec x_i^\mathrm{l}}\cdot\vec x^{\mathrm{l}}_{i, \vec y_{\vec{k}}} \, .
\end{equation}

%
\section{Examples \label{sec:examples}}
%
We proceed to demonstrate the application and efficacy of the proposed optimisation approach with four selected examples. In all the examples the compliance of the structure is the cost function and the material volume is prescribed. The first two examples are motivated by benchmark examples from continuum topology optimisation and concern the optimisation of lattice structures. In these two examples we consider only a lattice structure, without a shell skin, to study the proposed approach for lattice topology optimisation.
In the third example, we optimise the lattice infill topology of a lattice-skin cantilever structure.  The topology and shape optimisation of a lattice-skin  roof structure is considered in the last example. We use in all examples the sequential quadratic programming (SQP) optimisation algorithm in the NLopt library~\cite{johnsonnlopt}. The original cross-sectional areas of the struts are chosen as upper limits in lattice topology optimisation.

%
\subsection{Lattice cantilever} \label{sec:cantileverLattice}
\begin{figure}
\centering
  \includegraphics{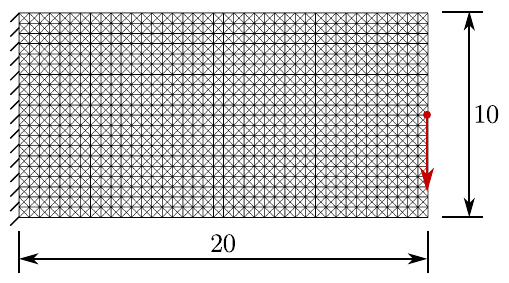}
	\caption{Lattice cantilever. Geometry, boundary conditions and loading. \label{fig:cantileverMiddle}}
\end{figure}

%
A 2D lattice cantilever of size $20 \times 10$ is considered, see Figure~\ref{fig:cantileverMiddle}. The left end of the lattice is fixed, while the remaining faces are free. A point load with magnitude 100 is applied at the mid-height on the right end. The periodic uniform lattice consists of square unit cells with two diagonals. That is, the four corners of the unit cell are connected to a centre joint by four struts. The side length of the unit cells  is~$0.5$, and the Young's modulus of the material is $E = 7 \times 10^7$. The total material volume of the original lattice is~$10$ with each strut having the same cross-sectional area~$5.1 \times 10^{-3}$. A volume fraction of $V_f^\mathrm{l} = 0.4$ is prescribed for optimisation. To this end, prior to  optimisation all cross-sectional areas are uniformly reduced to obtain a total material volume of~$4$. The cross-sectional areas of the struts are constrained not to exceed their initial value~$5.1 \times 10^{-3}$. The radius of the filter is~$R=1$, i.e. twice the unit cell side length.  To demonstrate the utility of the proposed SIMP-like penalisation approach, we consider two different types of functions for the penalisation of the relative strut densities~$\rho_e = A_e / \overline A_e$.

%
\subsubsection{Power function penalisation} \label{sec:cantileverLatticePower}
\begin{figure}
\centering
\subfloat[Penalisation factor~$p \in \{2, \, 3, \,4\}$, left to right. The respective compliances are~$J \in \{0.6019, \, 0.6306, \, 0.6404\}$.] 
{
  \centering
  \includegraphics{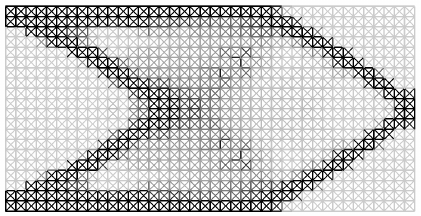} \hspace{2mm}
  \includegraphics{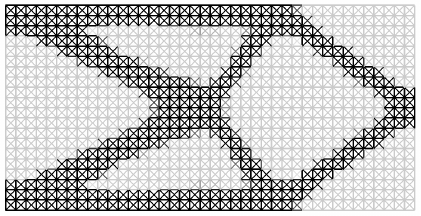} \hspace{2mm}
  \includegraphics{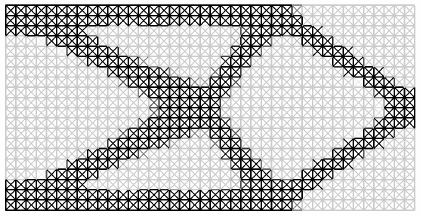}
  \label{fig:cantileverMiddlePenalty}
}
\\	
\subfloat[Frequency of cross-sectional areas for  $p \in \{2, \, 3, \,4\}$, left to right. The respective percentages of struts with intermediate areas are $48.4\%$, $14.4\%$ and $9.7\%$.] 
{
  \includegraphics[width = 0.3\textwidth]{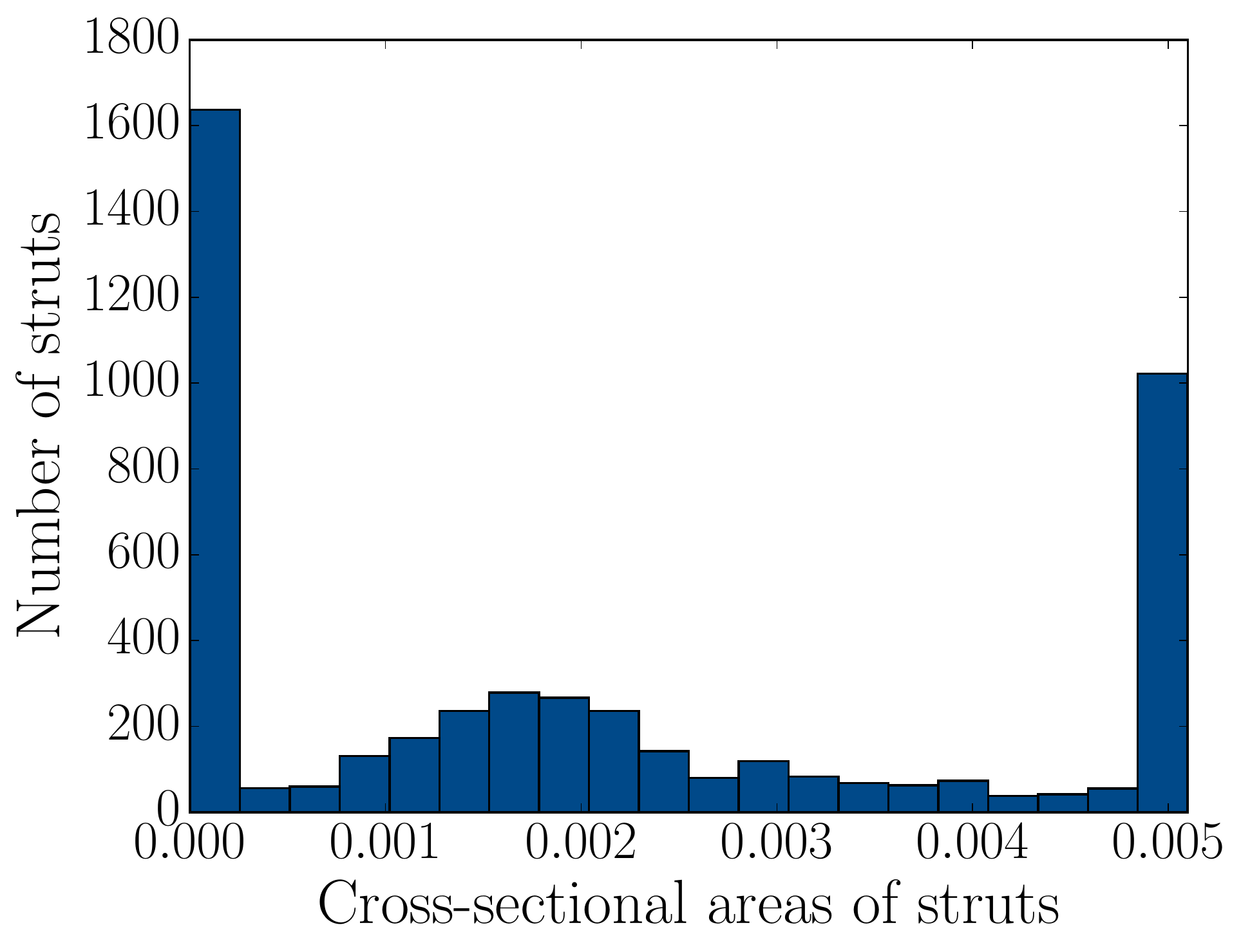}
  \hspace{2mm}
  \includegraphics[width = 0.3\textwidth]{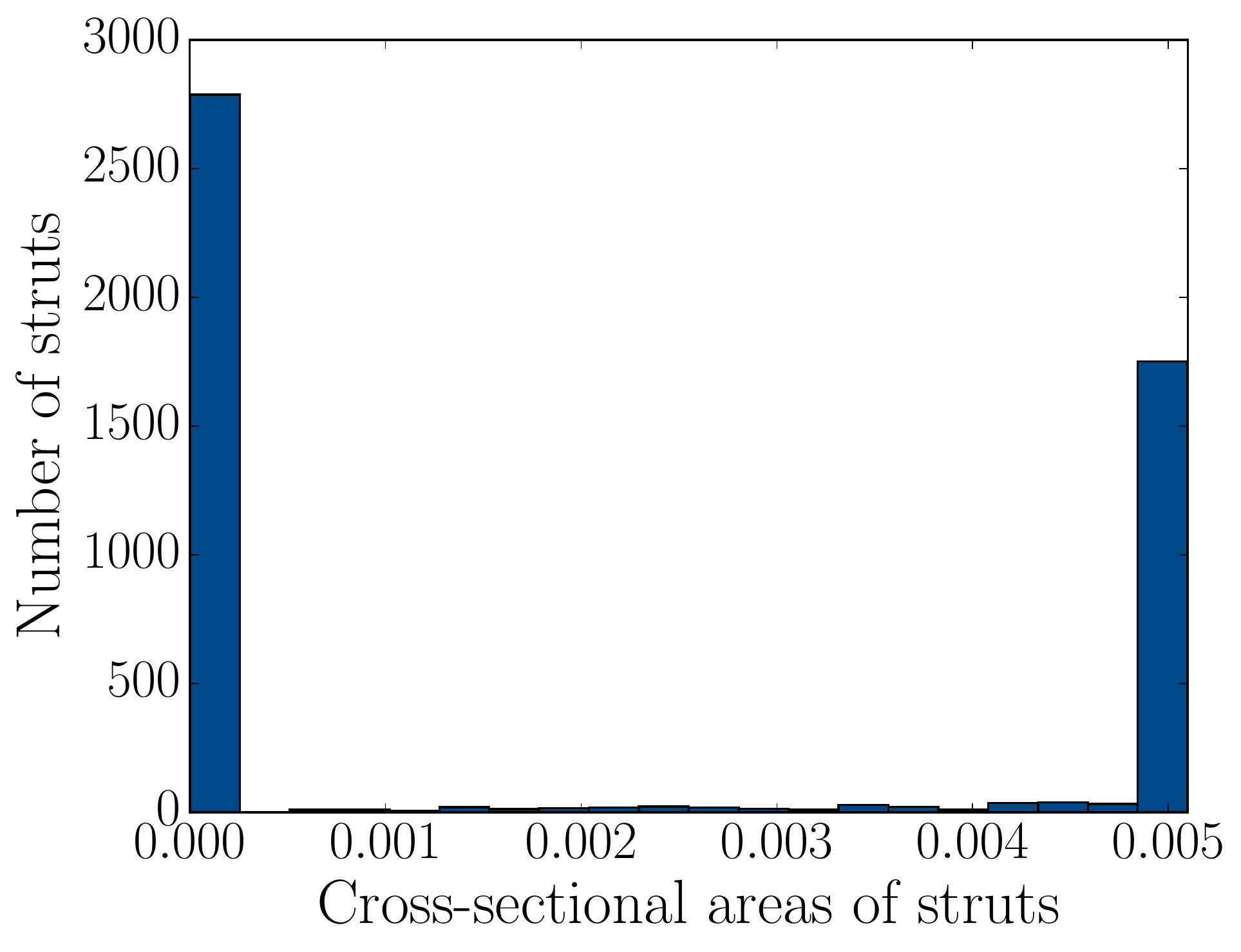}
  \hspace{2mm}
  \includegraphics[width = 0.3\textwidth]{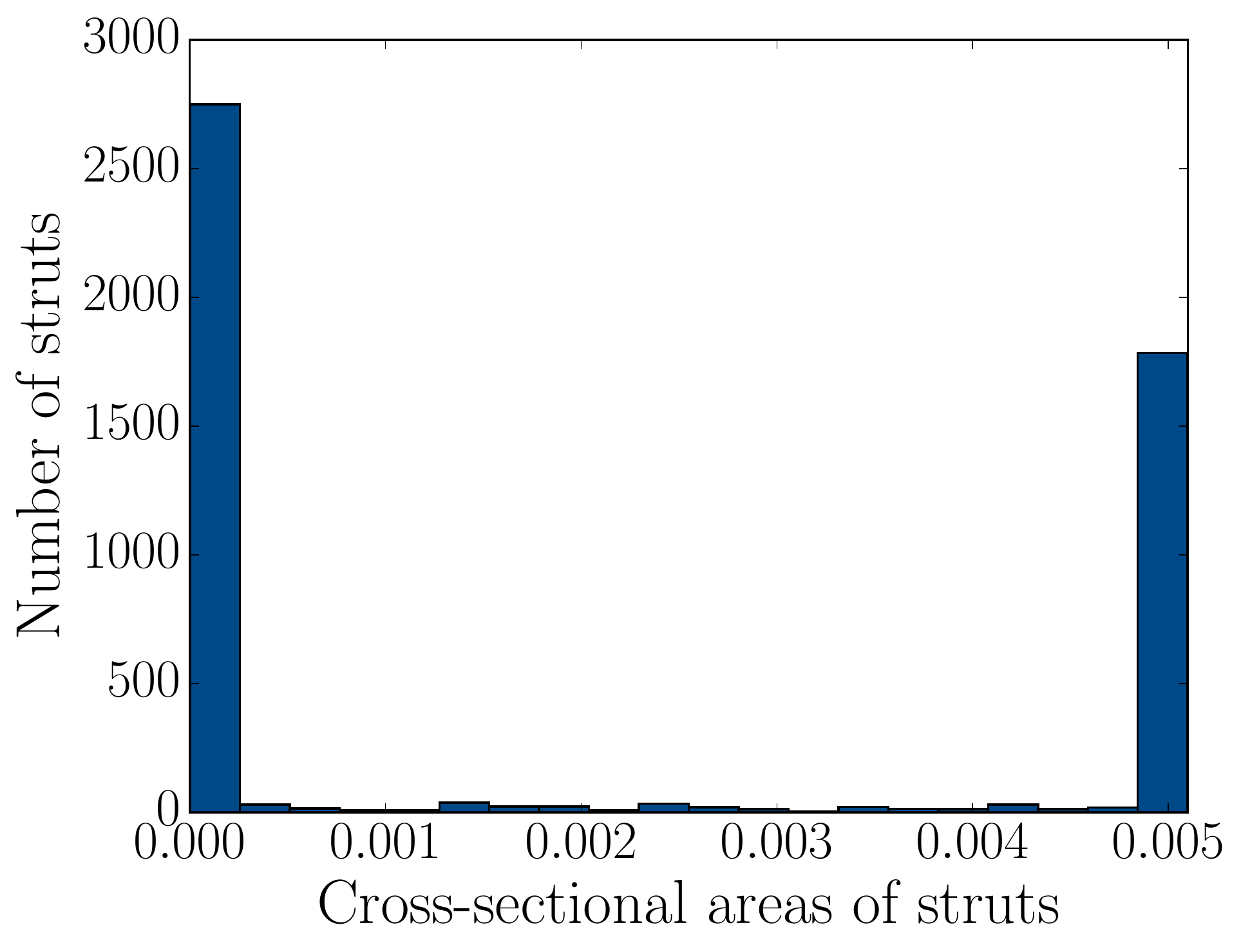}
  \label{fig:cantileverMiddlePenaltyAreas}
}
\caption{Topology optimisation of a lattice cantilever using power function penalisation and three different penalisation factors $p$.}
\end{figure}

%
First, the relative strut densities~$\rho_e$  are penalised with the power function introduced in~\eqref{eq:scaling}. We consider three different  penalisation parameters \mbox{$p \in \{ 2, \,  3, \, 4 \}$}, in turn, to investigate the influence of the choice of~$p$ on the results. The optimised cross-sectional areas of the struts are depicted in Figure~\ref{fig:cantileverMiddlePenalty}, and the frequency of the cross-sectional areas in Figure~\ref{fig:cantileverMiddlePenaltyAreas}. As can be seen, a too small penalisation factor leads to too many struts with intermediate cross-sectional areas. In addition, a too small penalisation factor may lead to checkerboarding as known from continuum optimisation.  Clearly, these issues can be alleviated by choosing a relatively large penalisation factor. Although the optimised structures with a lower penalisation achieve a smaller compliance, as given in the caption of Figure~\ref{fig:cantileverMiddlePenalty},  struts with deficiently small cross-sectional areas may buckle and can be impossible to manufacture.

\begin{figure}
\centering
\captionsetup[subfigure]{justification=centering}
\subfloat[Unit cell size $0.5 \times 0.5$; \newline $J = 0.6404$.]
{
  \includegraphics{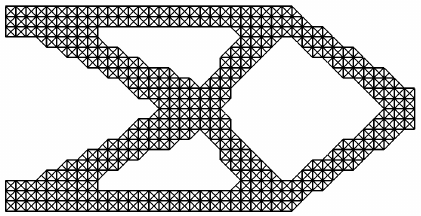}
}
\hfill
\subfloat[Unit cell size $0.25 \times 0.25$; \newline $J = 0.6506$.]
{
  \includegraphics{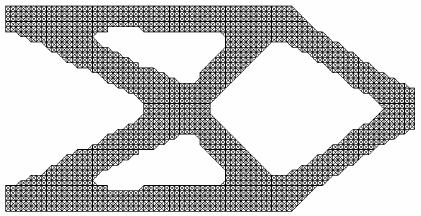}
}
\hfill
\subfloat[Unit cell size $0.125 \times 0.125$; \newline $J = 0.6643$.]
{
  \includegraphics{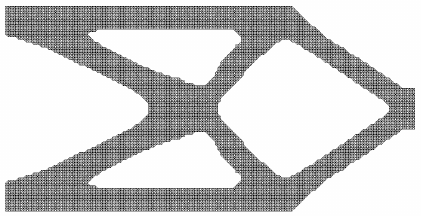}
}
\caption{Topology optimisation of a lattice cantilever using power function penalisation and three different unit cell sizes.}
\label{fig:cantileverMiddleDensity}
\end{figure}

We consider next three different lattices with the same total material volume, but with three different unit cell sizes \mbox{$0.5 \times 0.5$}, \mbox{$0.25 \times 0.25$} and \mbox{$ 0.125 \times 0.125$}. The penalisation factor is now fixed to $p = 4$. The optimised lattices shown in Figure~\ref{fig:cantileverMiddleDensity} are obtained by removing struts with cross-sectional areas less than~$0.001$ and recovering the complete topology of unit cells with any remaining diagonal struts, cf. Figure~\ref{fig:latticeExtraction}. As visually apparent the optimised topology and geometry remain the same irrespective of the unit cell size. However, as can be analytically explained the compliance becomes larger with decreasing unit cell size~\cite{pecullan1999scale}.

%
\subsubsection{B\'ezier penalisation} \label{sec:cantileverLatticeBezier}
\begin{figure}
	\centering
	\subfloat[Penalisation function and spline control  polygon (dashed) and control points~(green).  \label{fig:simpBezier}]{
	 \includegraphics[scale=1.35]{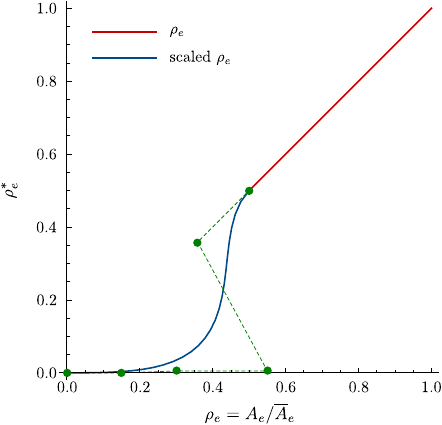}
}
	\hspace{0.1\textwidth}
	\subfloat[Three different penalisation functions.   \label{fig:cantileverPenaliseCurve}]
	{
	  \includegraphics[scale=1.3]{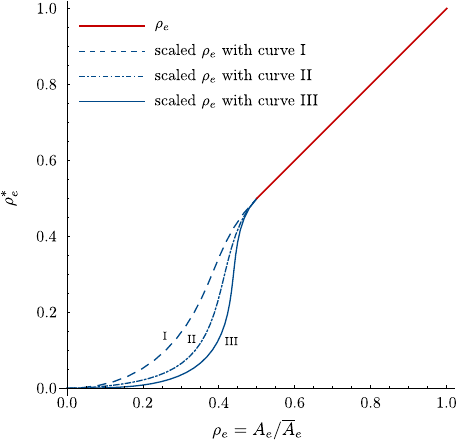}
	}
	\caption{Interpolation of the penalised density~$\rho_e^*$ with a B\'ezier spline curve (blue) and a smoothly connected line (red) at~$\rho_e=0.5$. \label{fig:BezierPenalisation} }
\end{figure}

\begin{figure}
\centering
\subfloat[Distribution of optimised cross-sectional areas for penalisation functions I, II and III, left to right, with respective compliances $J \in \{0.5683, \, 0.5851, \, 0.5933\}$.] 
{
  \centering
  \includegraphics[width=0.3\textwidth]{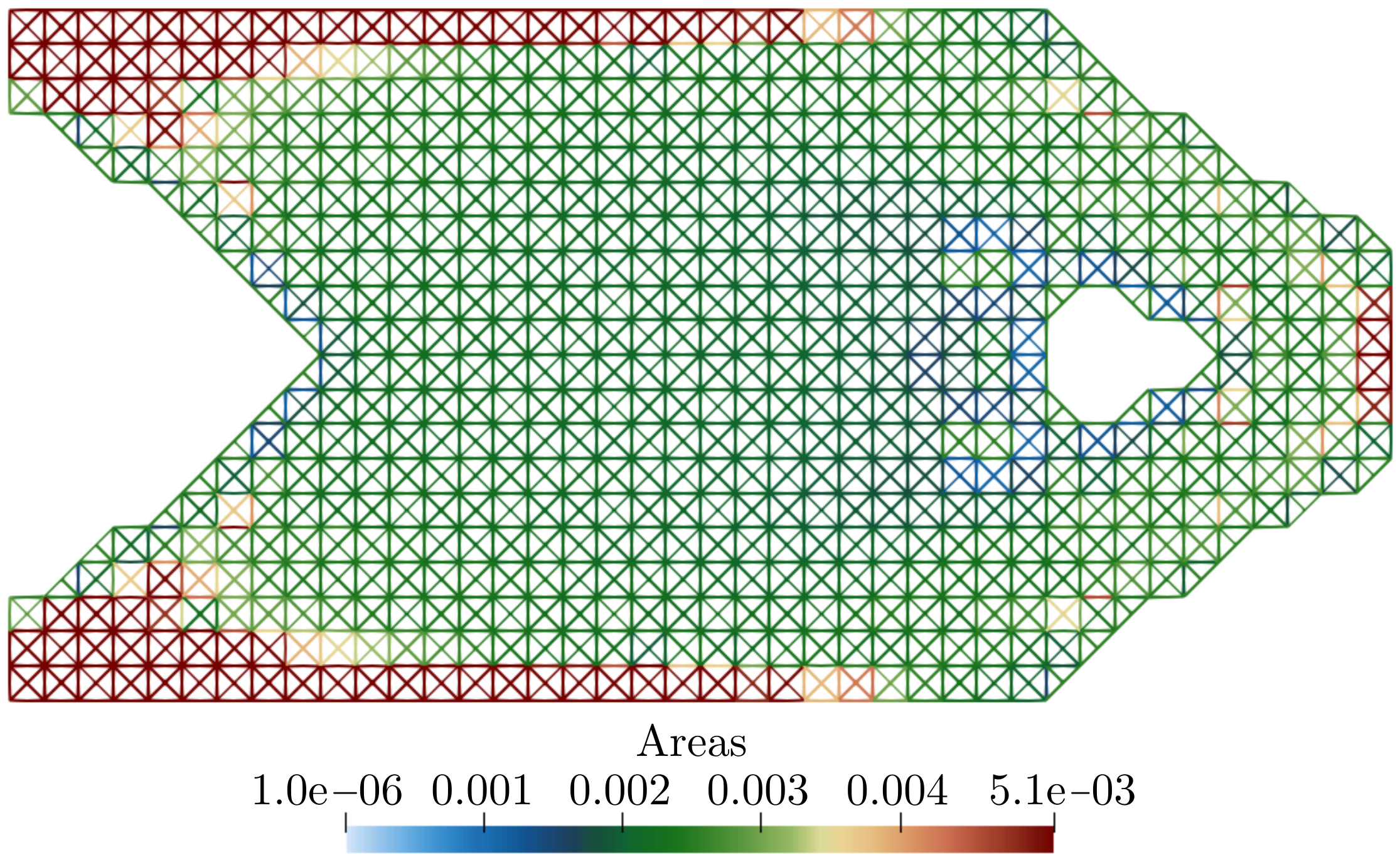} \hspace{4mm}
  \includegraphics[width=0.3\textwidth]{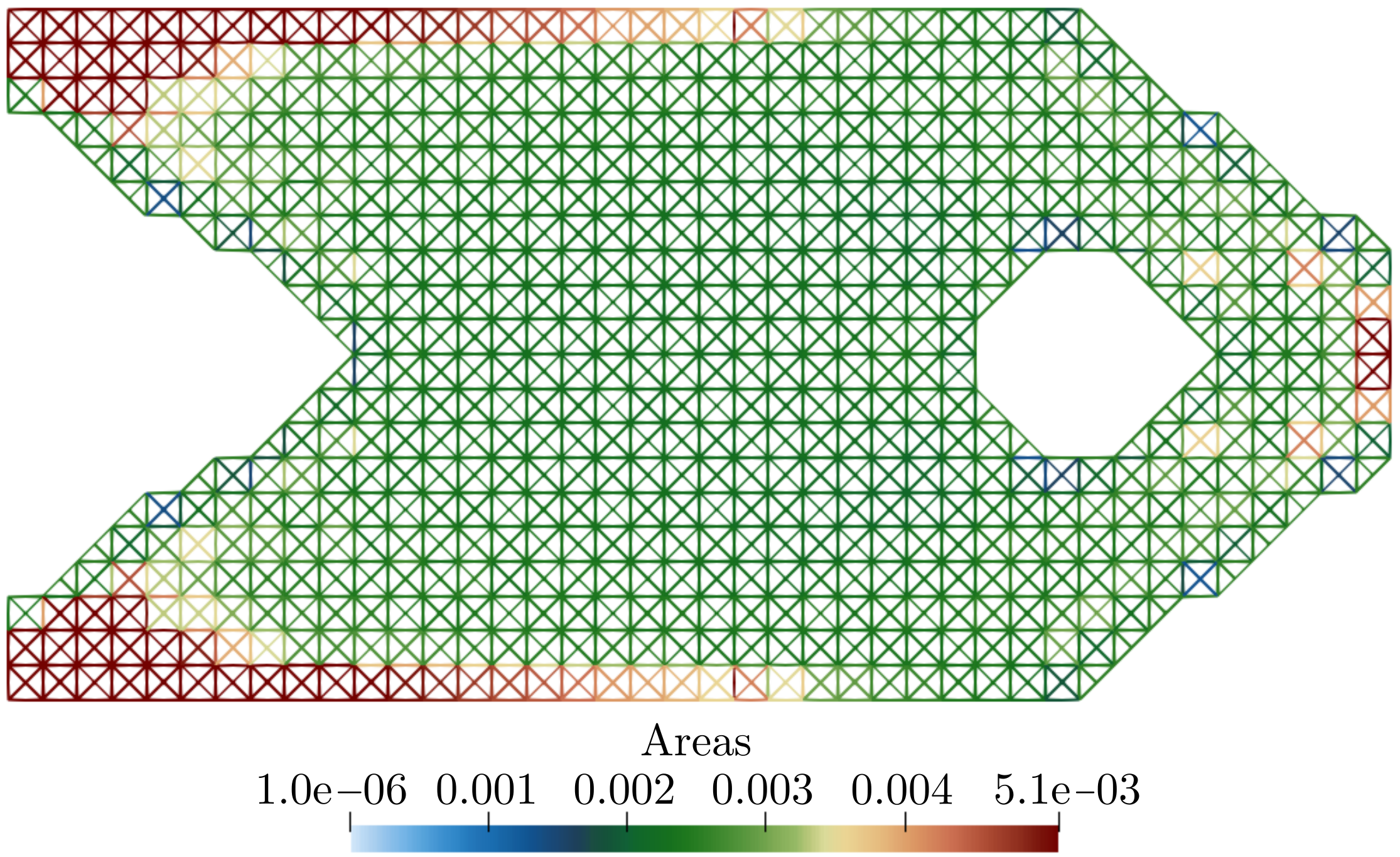} \hspace{4mm}
  \includegraphics[width=0.3\textwidth]{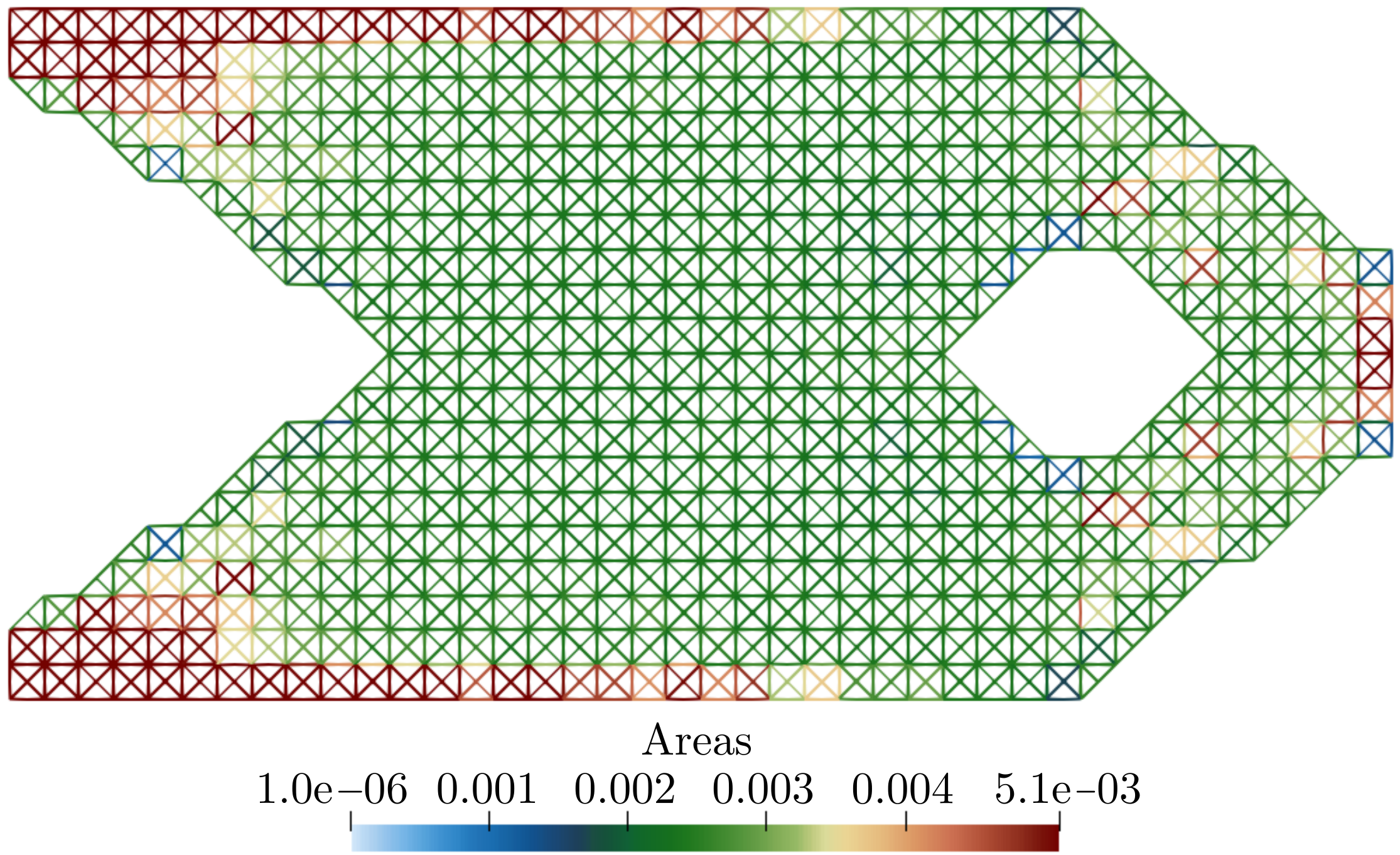}
  \label{fig:cantileverGraded}
}
\\	
\subfloat[Frequency of cross-sectional areas for penalisation functions I, II and III,  left to right. The respective percentage of intermediate areas with~$ 0 < A_e < 0.00255$ are $45.9\%$, $28.9\%$ and $9.8\%$.]
{
  \includegraphics[width=0.3\textwidth]{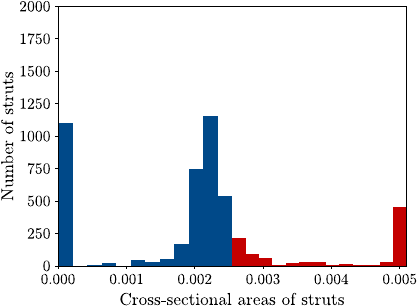} \hspace{4mm}
  \includegraphics[width=0.3\textwidth]{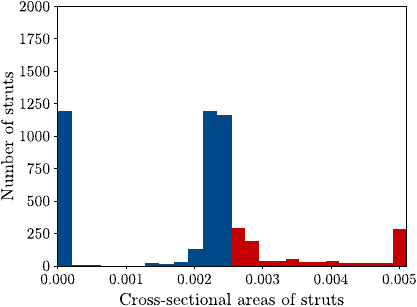} \hspace{4mm}
  \includegraphics[width=0.3\textwidth]{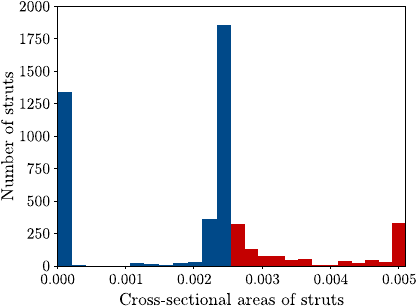}
  \label{fig:cantileverGradedHist}
}
\caption{Topology optimisation of a lattice cantilever using B\'ezier penalisation.}
\end{figure}

%
We choose as a penalisation function a B\'ezier spline of degree 5, which is smoothly connected to a line at~$\rho_e =0.5$, see Figure~\ref{fig:BezierPenalisation}. That is, only relative densities~$\rho_e \le 0.5$ are penalised ensuring that struts with relative densities close to~\mbox{$\rho_e=0$} and~\mbox{$\rho_e > 0.5$} are preferred. The relative density~$\rho_e = 0.5$ may be interpreted as the resolution lengthscale of the manufacturing process. We consider the spline curves I, II and III depicted in Figure~\ref{fig:cantileverPenaliseCurve}, in turn,  to investigate the influence of the shape of the spline curve on the optimisation results. The obtained optimised cross-sectional strut areas and their frequencies are given in Figure~\ref{fig:cantileverGraded} and~\ref{fig:cantileverGradedHist}, respectively. Struts with cross-sectional areas less than $0.001$ are removed from the lattice in Figure~\ref{fig:cantileverGraded}. Evidently, the grid-like lattice topology from before disappears due to the existence of intermediate cross-sectional areas, and the percentage of penalised cross-sectional areas~$ 0 < \rho_e < 0.5$  becomes smaller when a steeper function is used. 

%
\subsection{Lattice MBB beam}
\begin{figure}
\centering
\captionsetup[subfigure]{justification=centering}
\subfloat[Geometry, boundary conditions and loading.]
{
  \includegraphics[scale=0.9]{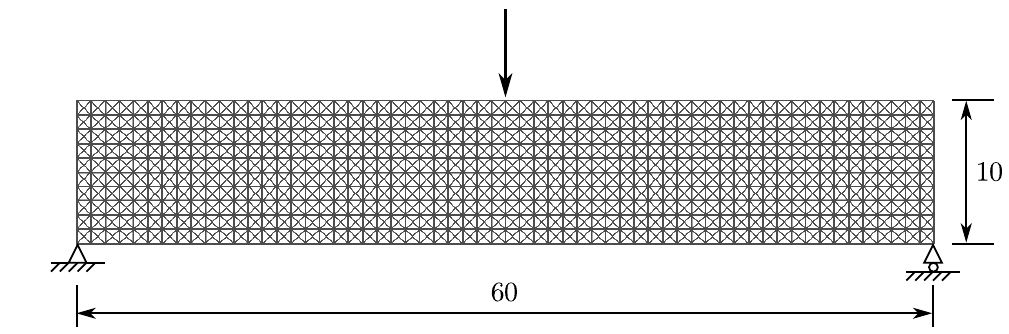}
  \label{fig:mbbInit}
}
\\
\subfloat[Cross-sectional areas of the struts after optimisation ($V_f^\mathrm{l} = 0.45$).] {
  \includegraphics[scale=0.9]{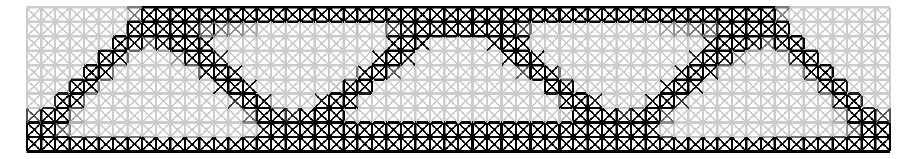}
  \label{fig:mbb1Opt}
}
\\
\subfloat[Final optimised structure ($J= 0.4689$).]
{
  \includegraphics[scale=0.9]{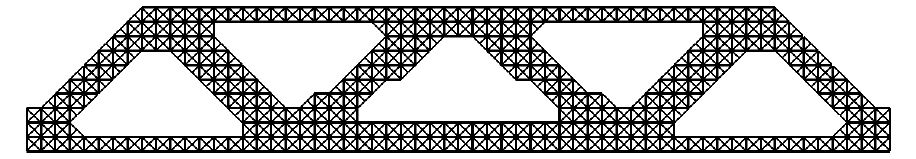}
  \label{fig:mbb1Final}
}
\caption{Topology optimisation of the lattice MBB structure.}
\end{figure}

%
The simply supported MBB beam is a widely studied benchmark example in continuum topology optimisation. We consider the equivalent lattice structure depicted in~Figure~\ref{fig:mbbInit}. The beam has the size $60 \times 10$ and is subjected to a vertical point load of $100$ at the midpoint of its top face. To begin with, we choose a periodic uniform lattice with unit cells of size $1 \times 1$ with each having two diagonals. The Young's modulus of the material is \mbox{$E = 7 \times 10^7$} and all the struts have the same cross-sectional area~$0.01685$ giving the total material volume~$50$. The compliance of the initial structure is $0.6563$. 

A volume fraction $V_f^\mathrm{l} = 0.45$ is prescribed for optimisation. Furthermore, using the power function penalisation the penalisation factor is chosen as $p = 3$ and the filter radius as $R=2$ (i.e. twice the unit cell size). The optimised cross-sectional areas of the struts are depicted in Figure~\ref{fig:mbb1Opt}. The final  structure in Figure~\ref{fig:mbb1Final} is obtained by removing struts with cross-sectional areas less than~$0.01$ and recovering the complete topology of unit cells which have remaining diagonal struts. The structural compliance of the optimised lattice is~$0.4689$, which represents a~$28.6\%$ reduction.

\begin{figure}
\centering
\captionsetup[subfigure]{justification=centering}
\subfloat[Geometry, boundary conditions and loading.]
{
  \hspace{0.3cm}
  \includegraphics[scale=0.9]{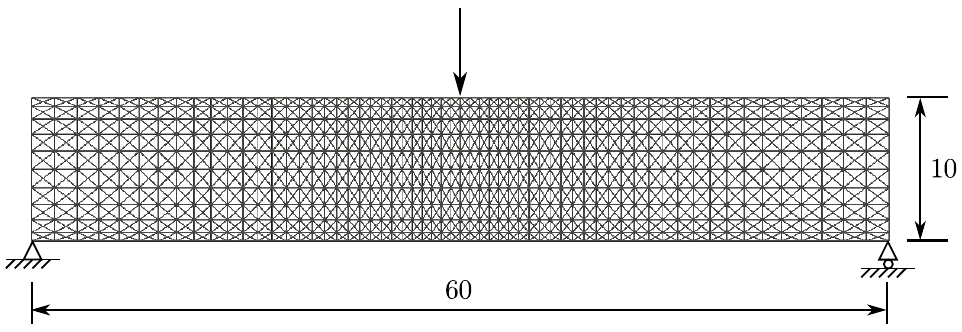}
  \label{fig:mbbNonUniformInit}
}
\\
\subfloat[Cross-sectional areas of the struts after optimisation ($V_f^\mathrm{l} = 0.45$).] 
{
  \includegraphics[scale=0.9]{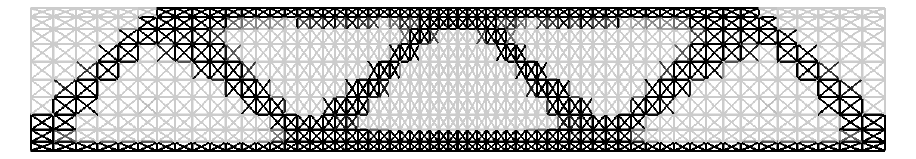}
  \label{fig:mbbNonUniformOpt}
}
\\
\subfloat[Final topology optimised structure ($J = 0.4349$).]
{
  \includegraphics[scale=0.9]{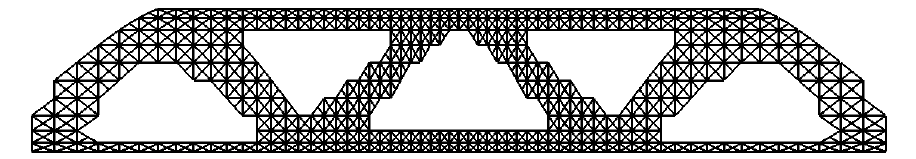}
  \label{fig:mbbNonUniformFinal}
}
\caption{Topology optimisation of the lattice MBB structure with a non-uniform distribution of unit cell sizes.}
\end{figure}

To demonstrate the application of the proposed approach to periodic non-uniform lattices, we consider the structure with non-uniform unit cell sizes shown in Figure~\ref{fig:mbbNonUniformInit}. The cells are designed to be densest at the mid-span where the axial forces are the largest. Again, a volume fraction of $45\%$ is chosen for  optimisation. The structural compliance of the initial lattice is~$0.5617$, which is reduced to~$0.4349$ by optimisation, see Figure~\ref{fig:mbbNonUniformFinal}. Compared to the uniform structure, the non-uniform structure is stiffer due to the purposefully chosen distribution of unit cell sizes, pointing out the potential of layout (shape) optimisation to further reduce compliance.

%
\subsection{Lattice-infilled cantilever}
\begin{figure}
\begin{minipage}{\textwidth}
\captionsetup[subfigure]{justification=centering}
\centering
\subfloat[Geometry, boundary conditions and loading.]
{
  \includegraphics[scale=0.8]{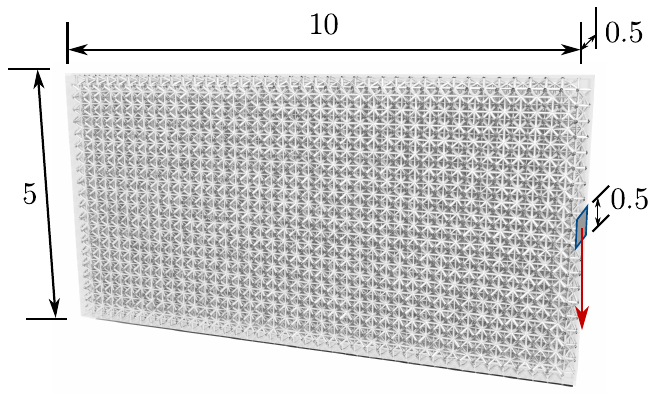}
  \label{fig:cantilever}
}
\hspace{1.5cm}
\subfloat[$V^\mathrm{l}/V^{\mathrm s} \approx 328.76$, $J = 0.3068$.]
{
  \includegraphics[scale=0.8]{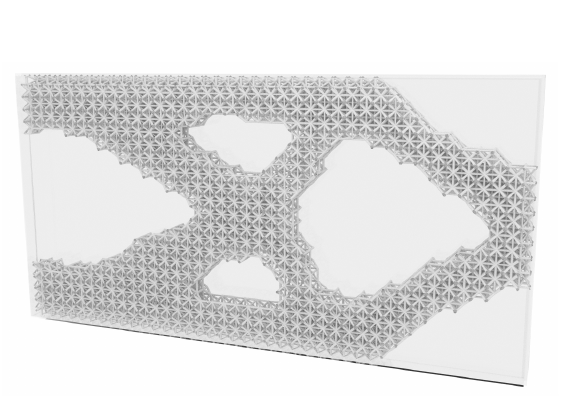}
  \label{fig:cantilever1}
}
\end{minipage}
\begin{minipage}{\textwidth}
\captionsetup[subfigure]{justification=centering}
\centering
\subfloat[$V^\mathrm{l}/V^{\mathrm s} \approx 20.99$,  $J = 0.2665$.]
{
  \includegraphics[scale=0.8]{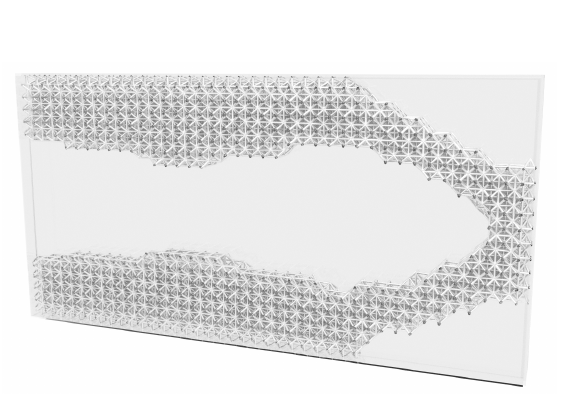}
  \label{fig:cantilever2}
}
\hspace{1.5cm}
\subfloat[$V^\mathrm{l}/V^{\mathrm s} \approx 2.3$, $J = 0.1364$.] 
{
  \includegraphics[scale=0.8]{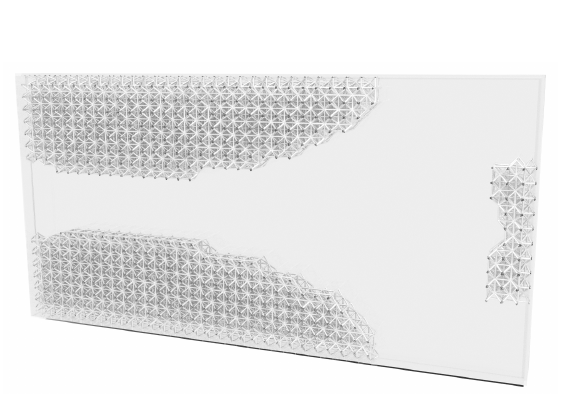}
  \label{fig:cantilever3}
}
\end{minipage}
\caption{Topology optimisation of a lattice-infilled cantilever considering three different lattice-to-shell volume ratios.}
\label{fig:cantileverOpt}
\end{figure}

%
We consider next the cantilever shown in Figure~\ref{fig:cantilever} comprised of a lattice and two exterior thin-shell face sheets. The left end is clamped while the remaining faces are free.  At the centre of the right face a distributed load with a total magnitude of~$200$ over an area of $0.5 \times 0.5$ is applied. The cantilever has the size \mbox{$10 \times 5 \times 0.5$}. The periodic uniform lattice consists of BCC (body-centred cubic) unit cells of side length~$0.25$. The Young's modulus is $E = 7 \times 10^7$ and the Poisson's ratio is~$0.35$. The total material volume is 40. The volume fraction is chosen as $V_f = 0.4$,  the penalisation factor in the power function penalisation is $p = 4$ and the filter radius as~$R=1$. In order to examine the effect of the thin-shell thickness on the structural performance, three different thicknesses~\mbox{$t \in \{ 0.001, \,  0.015, \, 0.1\}$} are considered. The total material volume is fixed so that the respective strut diameters are~$0.088$, $0.086$ and $0.074$. During the topology optimisation, the material volume ratio between the lattice and the thin-shell are preserved, i.e. the thickness of the thin-shells is fixed. 

The lattice topology optimisation results for the three different lattice-to-shell volume ratios $V^\mathrm{l}/V^{\mathrm s}$ are shown in Figures~\ref{fig:cantilever1},~\ref{fig:cantilever2} and~\ref{fig:cantilever3}. Clearly, the optimal lattice topology depends on~\mbox{$V^\mathrm{l}/V^{\mathrm s}$}. When the lattice volume is relatively large, the optimised structure resembles the ones obtained for lattice-only, c.f. Section~\ref{sec:cantileverLatticePower}, and continuum topology optimisation. However, with decreasing lattice volume the lattice-like topology starts to disappear, as can be seen in Figures~\ref{fig:cantilever2} and~\ref{fig:cantilever3}. This can be structurally explained as follows: when the shell thickness is very small, the shear rigidity in the height direction is provided by the lattice, whereas when the shell becomes thicker, the shear rigidity is increasingly provided by the shell face sheets. This leads to a concentration of  the struts at the top and bottom regions of the cantilever where the axial stresses due to bending are large.

It is worth emphasising that as the lattice volume decreases, the structural compliance of the optimised cantilever decreases as well, indicating that the lattice structure is not optimal compared with a solid-only structure in terms of compliance~\cite{sigmund2016non}. Nonetheless, a lattice can be optimal when in addition to compliance other design or performance criteria have to be taken into account, such as the buckling of the face sheets.

%
\subsection{Pentagon lattice-skin roof} 
\begin{figure}
\centering
\begin{minipage}{\textwidth}
\centering
\captionsetup[subfigure]{justification=centering}
\subfloat[Initial plate geometry with blue dots representing the supports.]
{
  \includegraphics[scale=0.068]{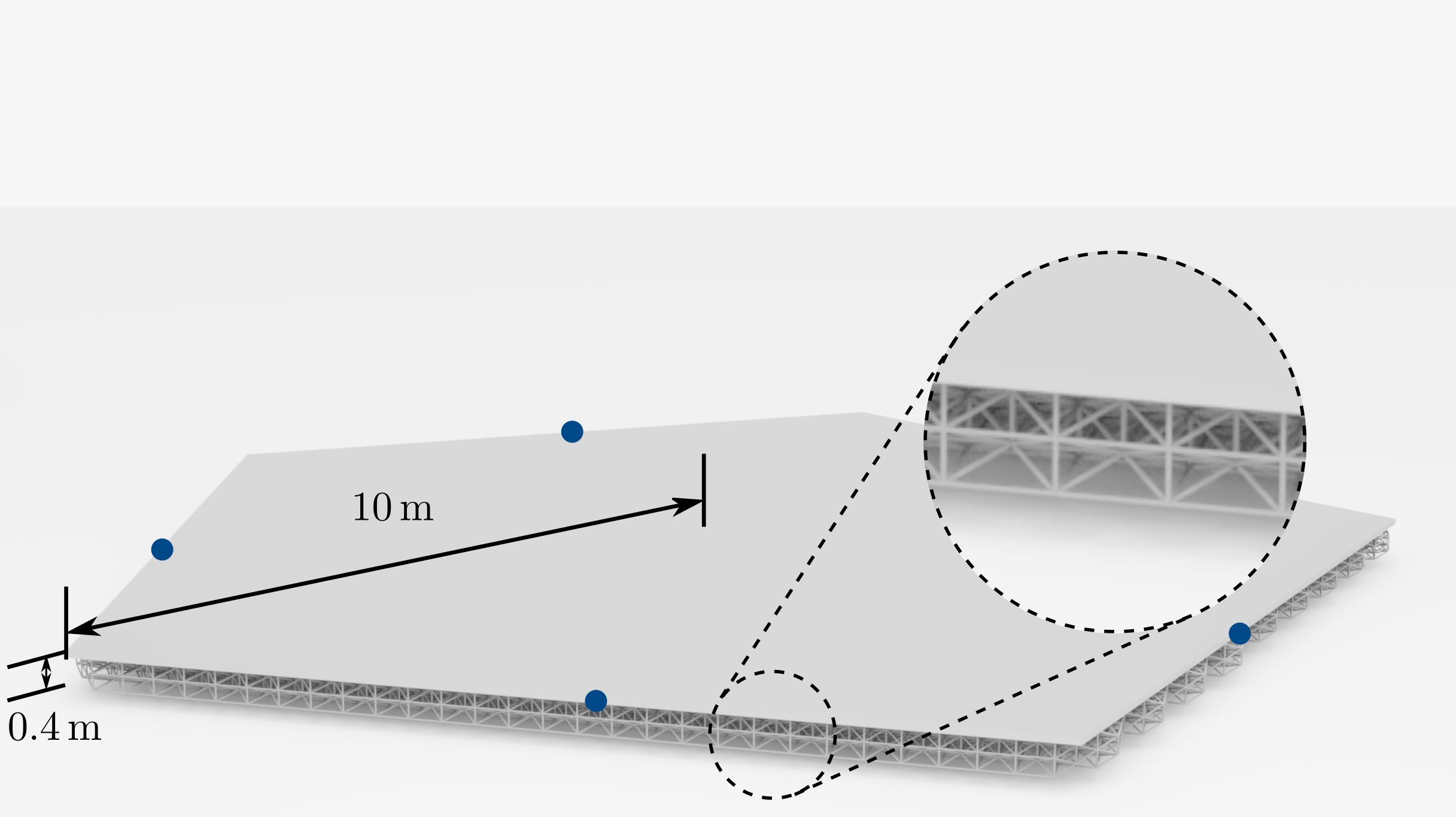}
  \label{fig:pentagonPlate}
}
\hspace{1cm}
\subfloat[Non-optimised curved roof after formfinding, $J = 159.44\,\mathrm{kN \, m}$.]
{
  \includegraphics{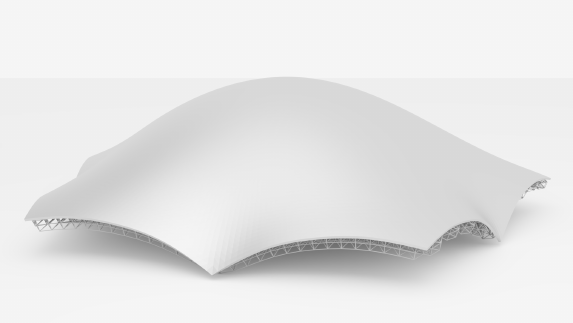}
  \label{fig:pentagonInit}
}
\vspace{0.5cm}
\end{minipage}
\begin{minipage}{\textwidth}
\centering
\captionsetup[subfigure]{justification=centering}
\subfloat[Shape optimised roof, $ J = 49.57\,\mathrm{kN \, m}$.] 
{
  \includegraphics{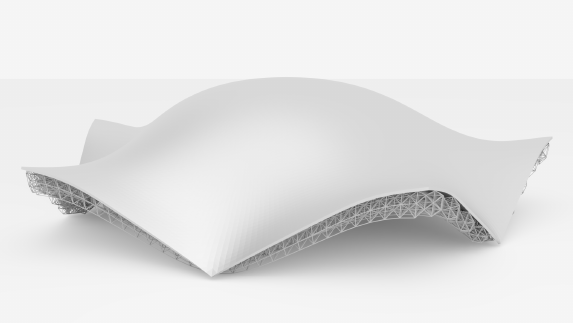}
  \label{fig:pentagonOpt}
}
\hspace{1cm}
\subfloat[Deflection of the shape optimised roof.]
{
  \includegraphics{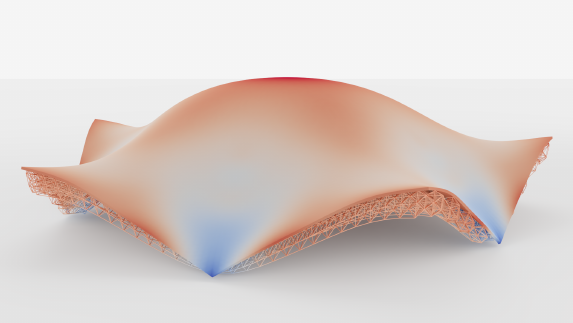}
  \label{fig:pentagonOptDisp}
}
\vspace{0.5cm}
\end{minipage}
\begin{minipage}{\textwidth}
\centering
\captionsetup[subfigure]{justification=centering}
\subfloat[A clipped view of the shape optimised roof.] 
{
  \includegraphics{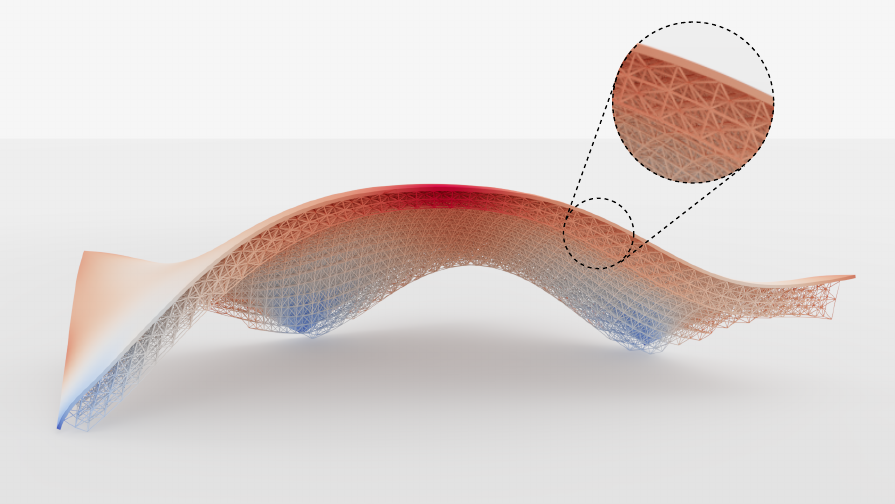}
  \label{fig:pentagonOptClip}
}
\end{minipage}
\caption{Shape optimisation of a lattice-skin roof subjected to vertical pressure loading.}
\end{figure}

%
As a last example we consider a lattice-skin roof to demonstrate the application of the proposed concurrent shape and lattice topology optimisation approach to more complex structures. The initial geometry of the curved roof is obtained starting from a pentagon-shaped flat shell, i.e. plate, supported by a lattice consisting of BCC periodic unit cells of size $0.4 \times 0.4 \times 0.2$, see Figure~\ref{fig:pentagonPlate}.  The radius of the circumcircle of the regular pentagon is $10\,\mathrm{m}$ and the shell thickness is $0.04\,\mathrm{m}$.  There are two unit cells through the thickness of the structure and in total $7493$ lattice joints and $39685$ struts. All the struts have a diameter of $0.01\,\mathrm{m}$. The Young's modulus and the Poisson's ratio are $E = 70\,\mathrm{GPa}$ and $\nu = 0.35$. A uniform pressure of~$50\,\mathrm{kN/m^2}$ is applied to the shell. The roof is supported  at the five pentagon midedges as indicated in Figure~\ref{fig:pentagonPlate} with blue solid dots.

Prior to optimisation, we obtain the curved roof in Figure~\ref{fig:pentagonInit} by scaling the deflection of the initially flat roof under a uniform pressure loading. This form-finding step is motivated by the hanging chain or cloth models commonly used in architectural design.  The so-obtained curved  lattice is guaranteed to be conformal to the curved shell mid-surface due to the coupling between the lattice and shell nodes. The compliance of the curved but not yet optimised roof is $159.44\,\mathrm{kN \, m}$. 

First, we optimise the shape of the roof after parameterising its geometry with the FFD technique presented in Section~\ref{sec:shapeOpt} and requiring that the total material volume is fixed. During iterative geometry updating, the positions of the supports are spatially fixed by placing FFD control points at the point supports  and constraining them to be fixed.  The final shape optimised roof is shown in Figure~\ref{fig:pentagonOpt}. The structural compliance is reduced to $49.57\,\mathrm{kN \, m}$, which gives a reduction of $68.9\%$ compared to the initial curved roof. Figure~\ref{fig:pentagonOptDisp} shows the displacement of the shape optimised roof under a uniform pressure loading. As can be seen in the clipped view in Figure~\ref{fig:pentagonOptClip}, in  the shape optimised roof the two-layer lattice is as required conformal to the thin-shell.

\begin{figure}
\centering
\begin{minipage}{\textwidth}
\centering
\captionsetup[subfigure]{justification=centering}
\subfloat[Initial shape optimised lattice,  $J = 53.90\,\mathrm{kN \, m}$.]
{
  \includegraphics[scale=0.9]{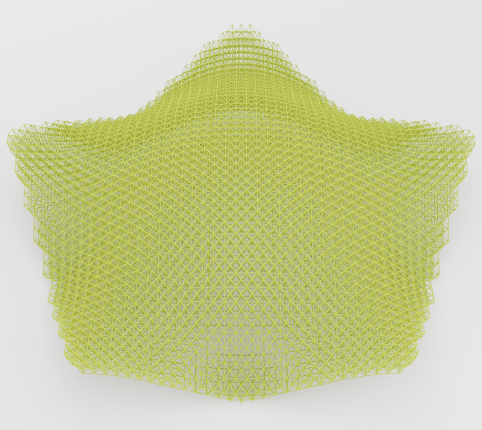}
  \label{fig:pentagonLattice}
}
\hfill
\subfloat[Topology and shape optimised lattice, $J = 51.24\, \mathrm{kN \, m}$.]
{
  \includegraphics[scale=0.9]{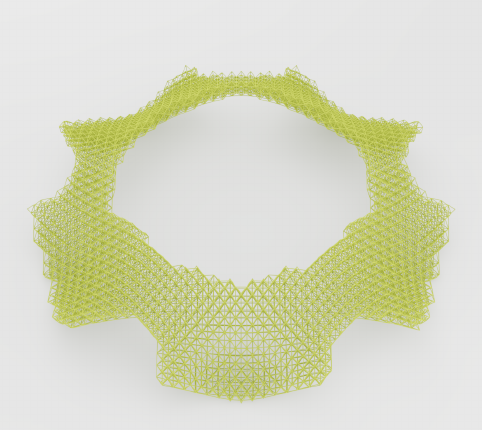}
  \label{fig:pentagonLatticeOpt}
}
\hfill
\subfloat[Deflection of the topology and shape optimised lattice.]
{
  \includegraphics[scale=0.9]{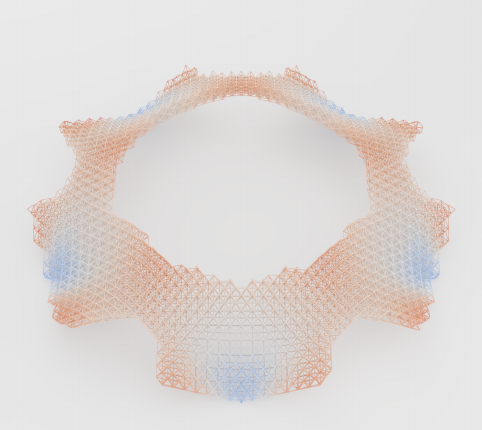}
  \label{fig:pentagonLatticeOptDisp}
}
\end{minipage}
\caption{Lattice topology optimisation of an already shape optimised lattice-skin roof. \label{fig:latticeTopoRoof}}
\end{figure}

\begin{figure}
\centering
  \includegraphics[scale=0.7]{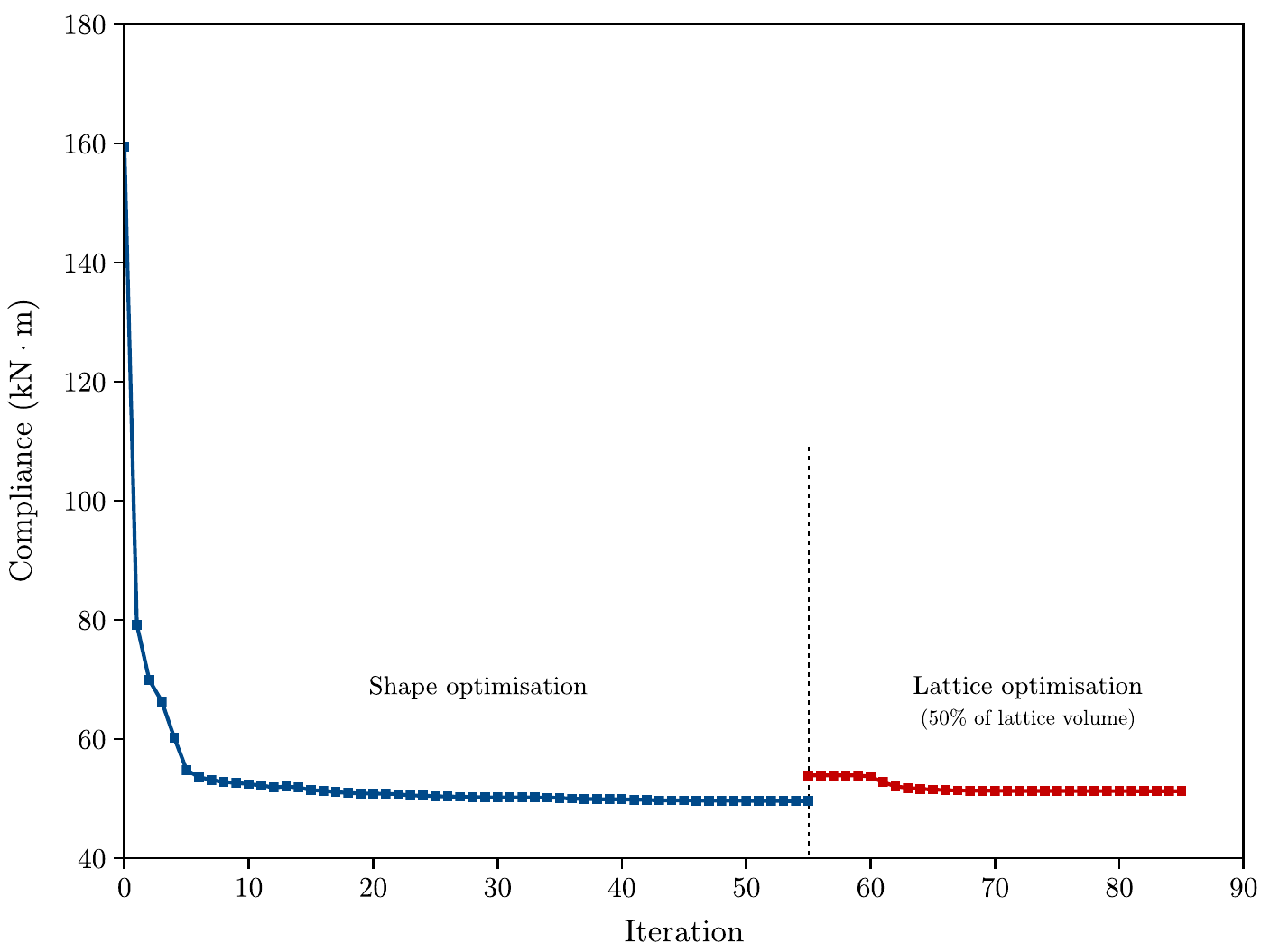}
\caption{Compliance of the lattice-skin roof during initial shape and subsequent lattice topology optimisation. \label{fig:latticeRoofConv}} 
\end{figure}

In a next step, the topology of the supporting lattice is optimised to reduce further the material usage without significantly compromising structural stiffness, see Figure~\ref{fig:latticeTopoRoof}. Note that the thin-shell is still taken into account in the respective sensitivity calculations. For lattice topology optimisation a volume fraction of $V^{\mathrm l}_f = 0.5$ is prescribed, where the cross-sectional areas of struts are set to be $50\%$ of the original values. The penalisation factor in the power function penalisation is $p = 4$ and the filter radius is $1\, \mathrm{m}$. To begin with, the cross-sectional areas of all struts are uniformly reduced without altering the lattice topology. As a consequence of the reduced lattice volume, the compliance jumps from $J=49.57\,\mathrm{kN \, m}$ to $53.9\,\mathrm{kN \, m}$.  The final shape and topology optimised roof is shown in Figure~\ref{fig:pentagonLatticeOpt}, which achieves a compliance of~$51.24\,\mathrm{kN \, m}$. The compliance increases slightly by $3.4\%$ compared with the shape optimised roof, while the lattice volume is reduced by $50\%$. Hence, the only shape optimised roof contains a large number of underutilised struts many of which can be discarded without compromising the stiffness of the roof greatly. As an ancillary benefit, the removal of struts may improve the natural lighting conditions in the space enclosed by the roof,  which is often an important architectural design consideration.  Figure~\ref{fig:pentagonLatticeOptDisp} shows the displacement of the final shape and topology optimised roof under uniform pressure loading. The decrease of the compliance during the initial shape and subsequent lattice topology optimisation are plotted in Figure~\ref{fig:latticeRoofConv}.

%
\section{Conclusions}
%

We introduced a novel shape and infill topology optimisation approach that fully takes into account the coupling between the pin-jointed lattice and the thin-shell skin. The modelling of the lattice as a pin-jointed truss allows us to relax the uniformity and scale separation assumptions underlying homogenisation theories. The proposed SIMP-like penalisation of cross-sectional areas yields struts with only desired areas, which is essential for satisfying manufacturability and other constraints, like strut buckling. In shape optimisation, we use the free-form deformation technique because it allows us to seamlessly parameterise the shape of the entire lattice-skin structure. While the shape and topology optimisation problems are solved in a sequential manner, the coupling between the lattice and the thin-shell is fully taken into account in both. The presented examples demonstrate the robustness and efficiency of the proposed approach in minimising the compliance of lattice and lattice-skin structures for a given volume of material. 

The presented lattice-skin optimisation approach can be extended in several ways. In this paper, we only optimised the compliance for a prescribed volume and a single load case. In practice, there are many more load cases, competing cost functions and constraints that have to be taken into account. We expect that most of these can be included without significant modifications to the presented approach~\cite{ohsaki2016optimization}. Furthermore, we modelled the lattice as a pin-jointed truss which is sufficient for engineered stretch-dominated lattices with appealing mechanical properties. Alternatively, to model bending-dominated lattices or to take into account secondary bending effects, the struts can be modelled as beams that are rigidly connected at the nodes~\cite{yin2020topologically}. This will increase the number of lattice degrees of freedom by a factor of two and will require the coupling of the rotations of the joints attached to the shell and the tangent plane rotations of the shell mid-surface. Finally, the manufacturability of the optimised geometries is a crucial consideration in practice~\cite{liu2016survey,liu2018current}. Although we did not consider any manufacturability constraints, a discrete lattice model should greatly aid their imposition.

%
\section*{Data Availability Statement}
The data that support the findings of this study are available from the corresponding author upon reasonable request.
%

\bibliographystyle{wileyj}
\bibliography{latticeSkinOpt}

\end{document}